\documentclass{amsart}
\usepackage[latin1]{inputenc}
\usepackage[english]{babel}
\usepackage{amsfonts}
\usepackage{amsmath}
\usepackage{amsthm,amssymb,latexsym}
\usepackage{graphicx}

\usepackage[active]{srcltx}

\usepackage{hyperref}

\newtheorem{theorem}{Theorem}[section]
\newtheorem{proposition}[theorem]{Proposition}
\newtheorem{lemma}[theorem]{Lemma}
\newtheorem{corollary}[theorem]{Corollary}

\newtheorem{remark}[theorem]{Remark}

\newcommand{\be}{\begin{enumerate}}
\newcommand{\ee}{\end{enumerate}}
\newcommand{\ba}{\begin{array}}
\newcommand{\ea}{\end{array}}
\newcommand{\bc}{\begin{center}}
\newcommand{\ec}{\end{center}}
\newcommand{\bq}{\begin{equation}}
\newcommand{\eq}{\end{equation}}

\newcommand{\F}{\mathbb F}
\newcommand{\A}{\mathbb A}
\newcommand{\Pp}{\mathbb P}
\newcommand{\K}{{\sf K}}

\newcommand{\fq}{\F_{\hskip-0.7mm q}}
\newcommand{\cfq}{\overline{\F}_{\hskip-0.7mm q}}
\newcommand{\fp}{\F_{\hskip-0.7mm p}}

\def\ifm#1#2{\relax \ifmmode#1\else#2\fi}

\begin{document}
\title[Symmetric hypersurfaces]{Singularities of symmetric hypersurfaces and an application to Reed-Solomon codes}%
\author{Antonio Cafure, Guillermo Matera, Melina Privitelli}%
\address{Instituto de Desarrollo Humano, UNGS}%
\email{acafure@ungs.edu.ar, gmatera@ungs.edu.ar, mprivitelli@conicet.gov.ar}%

\begin{abstract}
We determine conditions on $q$ for the nonexistence of deep holes of
the standard Reed--Solomon code of dimension $k$ over $\fq$
generated by polynomials of degree $k+d$. Our conditions rely on the
existence of $q$--rational points with nonzero, pairwise--distinct
coordinates of a certain family of hypersurfaces defined over $\fq$.
We show that the hypersurfaces under consideration are invariant
under the action of the symmetric group of permutations of the
coordinates. This allows us to obtain critical information
concerning the singular locus of these hypersurfaces, from which the
existence of $q$--rational points is established.
\end{abstract}
\maketitle
\section{Introduction}\label{section:introduction}
Let $\mathbb{F}_q$ be the finite field of $q$ elements of
characteristic $p$, let $\cfq$ denote its algebraic closure and let
$\fq^*$ denote the group of units of $\fq$. Let $\fq[T]$ and
$\fq[X_1,\ldots,X_n]$ denote the rings of univariate and $n$-variate
polynomials with coefficients in $\fq$, respectively.

Given  a subset  ${\sf D}:=\{x_1,\ldots,x_{n}\}\subset \fq$ and a
positive integer  $k\leq n$, the {\sf Reed--Solomon code  of length
$n$ and dimension $k$} over $\fq$ is the following subset of
$\fq^n$:
$$
{C}({\sf D},k):  = \{(f(x_{1}),\ldots,f(x_{n})) : f\in \fq[T],\,
\deg f \leq  k-1 \}.$$
The set ${\sf D}$ is called the {\sf evaluation set} and the
elements of ${C}({\sf D},k)$ are called {\sf codewords} of the code.
When  ${\sf D}= \fq^*$, we say that $C({\sf D}, k)$ is the {\sf
standard Reed--Solomon code}.

Let $C: = {C}({\sf D},k)$. For $\mathbf{w}\in \fq^{n}$,  we define
the {\sf distance of $\mathbf{w}$ to the code} $C$ as
$${\sf d}(\mathbf{w},C):=\displaystyle\min_{\mathbf{c}\in C}{\sf d}
(\mathbf{w},\mathbf{c})$$
where ${\sf d}$ is the Hamming distance of $\fq^n$. The {\sf minimum
distance} ${\sf d}(C)$ of $C$ is the shortest distance between any
two distinct codewords.  The {\sf covering radius} of $C$ is defined
as
$$
\rho:=\max_{\mathbf{y}\in\mathbb{F}_{\hskip-0.7mmq}^{n}}{\sf
d}(\mathbf{y},C).
$$
It is well--known that ${\sf d}(C)= n-k+1$ and $\rho= n-k$ holds.
Finally, we say that a word $\mathbf{w}\in \fq^{n}$ is a {\sf deep
hole} if ${\sf d}({\mathbf{w}},C)=\rho$ holds.

A decoding algorithm for the code $C$ receives a word
$\mathbf{w}\in\fq^{n}$ and outputs the message, namely the codeword
that is most likely to be received as $\mathbf{w}$ after
transmission, roughly speaking. One of the most important
algorithmic problems in this setting is that of the
\emph{maximum--likelihood decoding}, which consists in computing the
closest codeword to any given word $\mathbf{w}\in\fq^{n}$. It is
well--known that the maximum likelihood decoding problem for
Reed--Solomon codes is NP-complete (\cite{GuVa05}; see also
\cite{ChMu07}).

Suppose that we receive a word $\mathbf{w}:=(w_1,\ldots,w_n)\in
\fq^{n}$. Solving the maximum--likelihood decoding for $\mathbf{w}$
amounts at finding a polynomial $f \in \fq[T]$ of degree at most
$k-1$ satisfying the largest number of conditions $f(x_i)=w_i$,
$1\leq i\leq n$. By interpolation, there exists a unique polynomial
$f_\mathbf{w}$ of degree at most $n-1$ such that
$f_\mathbf{w}(x_i)=w_i$ holds for $1\leq i\leq n$. In this case, we
say that the word $\mathbf{w}$ was generated by the polynomial
$f_\mathbf{w}$. If $\deg f_\mathbf{w}\leq k-1$, then $\mathbf{w}$ is
a codeword.

In this paper, our main concern will be the existence of deep holes
of the given Reed--Solomon code $C$. According to our previous
remarks, a deep hole can only arise as the word generated by a
polynomial $f \in \fq[T]$ with $k\le \deg f\le n-1$. In this sense,
we have the following result.
\begin{proposition}\label{prop:poli de grado k genera deep holes}\cite[Corollary 1]{ChMu07}
Polynomials of degree $k$ generate deep holes.
\end{proposition}

Next we reduce further the set of polynomials $f\in\fq[T]$ which are
candidates for generating deep holes. Suppose that we receive a word
$\mathbf{w} \in \fq^n$, which is generated by a polynomial
$f_\mathbf{w}\in\fq[T]$ of degree greater than $k$. We want to know
whether $\mathbf{w}$ is a deep hole. We can decompose $f_\mathbf{w}$
as a sum $f_\mathbf{w}=g+h$, where $g$ consists of the sum of the
monomials of $f_\mathbf{w}$ of degree greater than or equal to $k$
and $h$ consists of those of degree less than or equal to $k-1$.
\begin{remark} \label{remark:observacion sobre las distancias}
Let $\mathbf{w}_g$ and $\mathbf{w}_h$ be the words generated by $g$
and $h$ respectively. Observe that $\mathbf{w}_h$ is a codeword. Let
$\mathbf{u}\in C$ be a codeword with ${\sf
d}(\mathbf{w},\mathbf{u})={\sf d}(\mathbf{w},C)$. From the
identities
$${\sf d}(\mathbf{w},C)={\sf d}(\mathbf{w},\mathbf{u})=
{\sf d}(\mathbf{w}-\mathbf{w}_h,\mathbf{u}-\mathbf{w}_h)= {\sf
d}(\mathbf{w}_g,\mathbf{u}-\mathbf{w}_h)$$
and the fact that  $\mathbf{u}-\mathbf{w}_h\in C$ holds, we conclude
$${\sf d}(\mathbf{w}_g,C)\le {\sf d}(\mathbf{w},C).$$
On the other hand, for $\mathbf{u}'\in C$ with ${\sf
d}(\mathbf{w}_g,C)={\sf d}(\mathbf{w}_g,\mathbf{u}')$, we have
$${\sf d}(\mathbf{w}_g,C) ={\sf d}(\mathbf{w}_g,\mathbf{u}') =
{\sf d}(\mathbf{w}_g+\mathbf{w}_h,\mathbf{u}' + \mathbf{w}_h)={\sf
d}(\mathbf{w},\mathbf{u}'+ \mathbf{w}_h)\ge {\sf d}(\mathbf{w},
C).$$
Therefore we have ${\sf d}(\mathbf{w},C) = {\sf d}(\mathbf{w}_g,C)$.
Hence $\mathbf{w}$ is a deep hole if and only if $\mathbf{w}_g$ is a
deep hole.
\end{remark}

From Remark \ref{remark:observacion sobre las distancias} it follows
that any deep hole of the Reed--Solomon code $C$ is generated by a
polynomial $f \in \fq[T]$ of the form
\begin{equation}\label{eq:polinomio f}
f:=T^{k+d}+f_{d-1}T^{k+d-1}+\cdots+f_{0}T^{k},
\end{equation}
where $d$ is a nonnegative integer with $k+d\le q-1$. In view of
Proposition \ref{prop:poli de grado k genera deep holes}, we shall
only discuss the case $d\ge 1$.

From now on we shall consider the standard Reed--Solomon code
$C:=C(\fq^*,k)$. In \cite{ChMu07} it is conjectured that the
reciprocal of Proposition \ref{prop:poli de grado k genera deep
holes} also holds, namely a word $\mathbf{w}$ is a deep hole of $C$
if and only if it is generated by a polynomial $f\in\fq[T]$ of
degree $k$. Furthermore, the existence of deep holes of $C$ is
related to the non-existence of $q$--rational points of a certain
family of hypersurfaces, in the way that we now explain. Fix
$f\in\fq[T]$ as in (\ref{eq:polinomio f})  and let  $\mathbf{w}_f$
be the generated word. Let $X_1,\ldots, X_{k+1}$ be indeterminates
over $\cfq$ and let $Q\in \fq[X_1,\ldots,X_{k+1}][T]$ be the
polynomial
$$Q=(T-X_1)\cdots(T-X_{k+1}).$$
We have that there exists $R_f\in\fq[X_1,\ldots,X_{k+1}][T]$ with
$\deg R_f\leq k$ such that the following relation holds:
\begin{equation}\label{eq:congruencia}
f\equiv R_f\mod{Q}.
\end{equation}
Assume that $R_f$ has degree $k$ and denote by $H_f\in
\fq[X_1,\ldots, X_{k+1}]$ its leading coefficient. Suppose that
there exists a vector $\mathbf{x}\in (\fq^*)^{k+1}$  with
pairwise--distinct coordinates such that $H_f(\mathbf{x})=0$ holds.
This implies that $r:=R_f(\mathbf{x},T)$ has degree at most $k-1$
and hence generates a codeword $\mathbf{w}_r$. By
(\ref{eq:congruencia}) we deduce that
$$d(\mathbf{w}_f,C) \leq d(\mathbf{w}_f,\mathbf{w}_r)\leq q-k-2$$
holds, and thus $\mathbf{w}_f$ is  not a deep hole.

As a consequence, we see that the given polynomial $f$ does not
generate a deep hole of $C$ if and only if there exists a zero
$\mathbf{x}:=(x_1,\ldots, x_{k+1})\in\fq^{k+1}$ of $H_f$ with nonzero,
pairwise--distinct coordinates, namely a solution
$\mathbf{x}\in\fq^{k+1}$ of the following system of equalities and
non-equalities:
\begin{equation}\label{eq:sistema de las soluciones no deseadas}
H_f (X_1,\ldots, X_{k+1})=0,\ \prod_{1\le i<j\le k+1}(X_i-X_j)\not=0,\
\prod_{1\le i\le k+1}X_i\not=0.
\end{equation}

\subsection{Related work}
As explained before, in \cite{ChMu07} the nonexistence of deep holes
of the standard Reed--Solomon code $C$ is reduced to the existence
of $q$--rational points, namely points whose coordinates belong to
$\fq$, with nonzero, pairwise--distinct coordinates of the
hypersurfaces $V_f$ defined by the family of polynomials $H_f$ of
(\ref{eq:sistema de las soluciones no deseadas}), where $f$ runs
through the set of polynomials $f\in\fq[T]$ as in (\ref{eq:polinomio
f}). The authors prove that all the hypersurfaces $V_f$ are
absolutely irreducible. This enables them to apply the explicit
version of the Lang--Weil estimate of \cite{CaMa06} in order to
obtain sufficient conditions for the nonexistence of deep holes of
Reed--Solomon codes. More precisely, the following result is
obtained.
\begin{theorem}
Let $k$, $d$ be given positive integers. If
$q>\max\{k^{7+\epsilon},d^{13/3+\epsilon}\}$ holds, then no word
$\mathbf{w}_f$ generated by a polynomial $f\in \fq[T]$ of degree
$k+d<q-1$ is a deep hole of the standard Reed--Solomon code over
$\fq$ of dimension $k$.
\end{theorem}

In \cite{LiWa08} the existence of deep holes is reconsidered. Using
the Weil estimate for certain character sums as in \cite{Wan97}, the
authors obtain the following result.
\begin{theorem}\label{teo:li-wan} Let $k$, $d$ be given positive
integers. If $q>\max\{d^{2+\epsilon},(k+1)^2\}$ and
$k>(\frac{2}{\epsilon}+1)d+ \frac{8}{\epsilon}+2$ holds for a
constant $\epsilon>0$, then no word $\mathbf{w}_f$ generated by a
polynomial $f\in \fq[T]$ of degree $k+d<q-1$ is a deep hole of the
standard Reed--Solomon code over $\fq$ of dimension $k$.
 \end{theorem}

\subsection{Our results}
We determine further threshold values $\lambda_1(d,k)$ and
$\lambda_2(d)$ such that for $q>\lambda_1(d,k)$ and $k>\lambda_2(d)$
the standard Reed--Solomon code over $\fq$ of dimension $k$ has no
deep holes generated by polynomials of degree $k+d$. In fact, we
have the following result (see Theorems \ref{teo:principal - version
precisa} and \ref{teo:principal char greater than d} for precise
versions).
\begin{theorem} \label{teorema principal}
Let $k$, $d$ be positive integers and $0<\epsilon<1$. Suppose that
$q>\{d^{2+\epsilon},(k+1)^2\}$ and $ k>(\frac{2}{\epsilon}+1)d$
hold. Let $f\in \fq[T]$ be an arbitrary polynomial of degree
$k+d<q-1$ and let $\mathbf{w}_f\in \fq^{q-1}$ be the word defined by
$f$. Then $\mathbf{w}_f$ is not a deep hole of the standard
Reed--Solomon code over $\fq$ of dimension $k$.
\end{theorem}

This result is obtained from a lower bound on the number of
$q$--rational points with nonzero, pairwise--distinct coordinates of
the family of hypersurfaces $V_f$ introduced above. Our result
improves that of \cite{ChMu07} by means of a deeper study of the
geometry of these hypersurfaces. In fact, we show that each
hypersurface $V_f$ has a singular locus of dimension at most $d-1$
(Corollary \ref{coro:dimension singular locus Vf}), which in
particular implies that it is absolutely irreducible (as proved by
\cite{ChMu07}). We further prove that for $\mathrm{char}\, \fq>d+1$,
the singular locus of the hypersurfaces $V_f$ of interest has
dimension at most $d-2$ (Theorem \ref{teo:singular locus dim d-1
implies monomial} and Proposition \ref{prop:monomial doesnt generate
deep holes}).

For this purpose, we show that the polynomials $H_f\in
\fq[X_1,\ldots,X_{k+1}]$ defining the hypersurfaces $V_f$ are
symmetric, namely invariant under any permutation of the variables
$X_1,\ldots, X_{k+1}$. More precisely, for any polynomial $f\in\fq[T]$
as in (\ref{eq:polinomio f}) of degree $d$, we prove that $H_f$ can
be expressed as a polynomial in the first $d$ elementary symmetric
polynomials $\Pi_1,\ldots,\Pi_d$ of $\fq[X_1,\ldots, X_{k+1}]$ (Propositions
\ref{Proposicion: formula explicita para la hipersuperficie H_0d}
and \ref{prop: formula para Hf}). Such an expression involves the
number of different partitions of $d$ (admitting repetition) and
resembles the Waring formula.

The result on the dimension of the singular locus of the
hypersurfaces $V_f$ is then combined with estimates on the number of
$q$--rational points of singular complete intersections
\cite{GhLa02a}, yielding our main result Theorem \ref{teorema
principal}.

Our results also constitute an improvement of that of \cite{LiWa08},
as can be readily deduced by comparing the statements of Theorems
\ref{teo:li-wan} and \ref{teorema principal}. Nevertheless, as the
``main'' exponents in both results are similar, we would like to
stress here the methodological aspect. As mentioned before, the
critical point for our approach is the invariance of the family of
hypersurfaces $V_f$ under the action of the symmetric group of $k+1$
elements. In fact, our results on the dimension of the singular
locus and the estimates on the number of $q$--rational points can be
extended \textit{mutatis mutandis} to any symmetric hypersurface
whose projection on the set of primary invariants (using the
terminology of invariant theory) defines a nonsingular hypersurface.
This might be seen as a further source of interest of our approach,
since hypersurfaces with symmetries arise frequently in coding
theory and cryptography (for example, in the study of almost perfect
nonlinear polynomials or differentially uniform mappings; see, e.g.,
\cite{Rodier09} or \cite{AuRo09}).

\section{$H_f$ in terms of the elementary symmetric
polynomials} \label{section:planteo del problema}

Fix positive integers $d$ and $k$ such that $d < k$, and consider
the first $d$ elementary symmetric polynomials $\Pi_1 ,\ldots,
\Pi_d$ of $\fq[X_1,\ldots, X_{k+1}]$. In Section
\ref{section:introduction} we associate a polynomial $H_f\in
\fq[X_1,\ldots, X_{k+1}]$ to every polynomial $f\in\fq[T]$ of degree
$k+d$ as in (\ref{eq:polinomio f}). As asserted above, the word
$\mathbf{w}_f$ generated by a given polynomial $f$ is not a deep
hole of the standard Reed--Solomon code of dimension $k$ over $\fq$
if $H_f$ has a $q$--rational zero with nonzero, pairwise--distinct
coordinates.

The main purpose of this section is to show how the polynomials
$H_f$ can be expressed in terms of the elementary symmetric
polynomials $\Pi_1 ,\ldots, \Pi_d$. For this purpose, we first obtain a
recursive expression for the polynomial $H_d$ associated to the
monomial $T^{k+d}$.

\begin{lemma}\label{lema: formula recursiva para los H0d}
Fix $H_0:=1$. For any $d \geq 1$, the following identity holds:

\begin{equation}
\label{eq:formula recursiva para H_0d}
H_d=\Pi_1H_{d-1}-\Pi_2H_{d-2}+\cdots+(-1)^{d-1}\Pi_dH_0.
\end{equation}
\end{lemma}
\begin{proof}
Let as before $Q:=(T-X_1)\cdots(T-X_{k+1})$. We have

$$
T^{k+1}\equiv  \Pi_1T^k-\Pi_2T^{k-1}\!+\cdots+(-1)^{d-1}\Pi_d
T^{k-(d-1)}+\cdots+(-1)^{k}\Pi_{k+1} \mod{Q}.$$

Multiplying this congruence relation by $T^{d-1}$ we obtain:
$$T^{k+d}\equiv \Pi_1T^{k+d-1}-\Pi_2T^{k+d-2}+\cdots+(-1)^{d-1}
\Pi_{d}T^k+\mathcal{O}(T^{k-1})\mod{Q},$$
where $\mathcal{O}(T^{k-1})$ represents a sum of terms of
$\fq[X_1,\ldots,X_{k+1}][T]$ of degree at most $k-1$ in $T$. Recall
that we define $H_{d-j}$ as the unique polynomial of
$\fq[X_1,\ldots, X_{k+1}]$ satisfying the congruence relation
$$T^{k+d-j}\equiv H_{d-j}T^k+\mathcal{O}(T^{k-1})\mod{Q}$$
for $1\leq j\leq d-1$. Hence, we obtain the equality
$$H_d=\Pi_1 H_{d-1}-\Pi_2H_{d-2}+\cdots+(-1)^{d-1}\Pi_d.$$
This finishes the proof of the lemma.
\end{proof}

Our second step is to obtain an explicit expression of the
polynomial $H_d$ in terms of the elementary symmetric polynomials
$\Pi_1,\ldots,\Pi_d$. From this expression we readily obtain an
expression for the polynomial $H_f$ associated to an arbitrary
polynomial $f$ as in (\ref{eq:polinomio f}) of degree $k+d$.
\begin{proposition}\label{Proposicion: formula explicita para la hipersuperficie H_0d}
Let $H_d\in\fq[X_1,\ldots, X_{k+1}]$ be the polynomial associated to the
monomial $T^{k+d}$. Then the following identity holds:
\begin{equation}\label{eq:formula explicita para los HOd}
H_d = \sum_{i_1+2i_2+\cdots+ di_d=d} (-1)^{\Delta(i_1,\ldots, i_d)}
\dfrac{(i_1+\cdots+ i_d)!}{i_1!\cdots i_d!}\Pi_1^{i_1}\cdots
          \Pi_d^{i_d},
\end{equation}
where $0\le i_j\le d$ holds for $1\le j\le d$ and $\Delta(i_1,\ldots,
i_d):=i_2+i_4+\cdots+ i_{2\lfloor d/2\rfloor}$ denotes the sum of
indices $i_j$ for which $j$ is an even number.
\end{proposition}

\begin{proof}
We argue by induction on $d$. The case $d=1$ follows immediately
from (\ref{eq:formula recursiva para H_0d}).

Assume now that $d>1$ holds and (\ref{eq:formula explicita para los
HOd}) is valid for $1\leq j \leq d-1$. From (\ref{eq:formula
explicita para los HOd}) we easily conclude that $H_j$ is a
homogeneous symmetric polynomial of $\fq[X_1,\ldots, X_{k+1}]$ of degree $j$
for $1\le j\le d-1$. Furthermore, from Lemma \ref{lema: formula
recursiva para los H0d} we deduce that $H_d$ is also a homogeneous
symmetric polynomial of degree $d$. Combining the inductive
hypotheses and Lemma \ref{lema: formula recursiva para los H0d} we
see that $H_d$ can be expressed in the form
$$H_d=\sum_{i_1+\cdots+ di_d=d}
a_{i_1,\ldots,i_d}\Pi_1^{i_1} \cdots \Pi_d^{i_d},$$
for suitable elements $a_{i_1,\ldots, i_d}\in\fq$. As a consequence, it
only remains to prove that the terms $a_{i_1,\ldots, i_d}$ have the
asserted form, namely
$$a_{i_1,\ldots, i_d}=(-1)^{\Delta(i_1,\ldots, i_d)}
\dfrac{(i_1+\cdots+ i_d)!}{i_1!\cdots i_d!}.$$

Fix $(i_1,\ldots,i_d)\in(\mathbb{Z}_{\ge 0})^d$ with $i_1+2i_2+\cdots+
d{i_d}=d$. Then Lemma \ref{lema: formula recursiva para los H0d}
shows that
$$a_{i_1,\ldots, i_d}=\sum_{j=1}^d(-1)^{j-1}
(H_{d-j})_{i_1,\ldots, i_j-1,\ldots, i_d},$$
where $(H_{d-j})_{{i_1},\ldots, {i_j-1},\ldots, {i_d}}$ is the coefficient
of the monomial $\Pi_1^{i_1}\cdots \Pi_j^{i_j-1}\cdots \Pi_d^{i_d}$
in the expression of $H_{d-j}$ as a polynomial of
$\fq[\Pi_1,\ldots,\Pi_d]$.

Therefore, applying the inductive hypothesis, we obtain:
\begin{eqnarray*}
a_{i_1,\ldots, i_d}= \sum_{j=1}^d (-1)^{j-1}(-1)^{\Delta(i_1,\ldots,
{i_j-1},\ldots, i_d)}\,\dfrac{(i_1+\cdots+ {i_d}-1)!}{i_1!\cdots
(i_j-1)!\cdots i_d!}.
\end{eqnarray*}
If $j$ is an odd number, then $\Delta(i_1,\ldots, {i_j}-1 ,\ldots,
i_d)=\Delta({i_1},\ldots, {i_j},\ldots,{i_d})$ and $(-1)^{j-1}=1$ hold,
which implies $(-1)^{j-1+\Delta(i_1,\ldots,{i_j}-1,\ldots,{i_d})}=
(-1)^{\Delta({i_1},\ldots,{i_j},\ldots,{i_d})}$. On the other hand, if $j$
is an even number then we have $(-1)^{{j}-1}=-1$ and
$(-1)^{\Delta({i_1} ,\ldots,
{i_j},\ldots,{i_d})}=(-1)^{{j}-1}(-1)^{\Delta({i_1} ,\ldots, {i_j}-1 ,\ldots,
{i_d})}$. Therefore
\begin{eqnarray*}
a_{{i_1},\ldots,{i_d}}& = &
(-1)^{\Delta({i_1},\ldots,{i_d})}({i_1}+\cdots+{i_d}-1)!
\dfrac{({i_1}+\cdots+{i_d})}{{i_1}!\ldots
          {i_d}!}\\
& = & (-1)^{\Delta({i_1},\ldots,{i_d})}
          \dfrac{({i_1}+\cdots +{i_d})!}{{i_1}!\ldots
          {i_d}!}.
\end{eqnarray*}
This concludes the proof of the proposition.
\end{proof}

It is interesting to remark the similarity of the expression for
$H_d$ with Waring's formula expressing the power sums in terms of
the elementary symmetric polynomials (see, e.g., \cite[Theorem
1.76]{LiNi97}).

Finally we obtain an expression of the polynomial $H_f\in
\fq[X_1,\ldots,X_{k+1}]$ associated to an arbitrary polynomial
$f\in\fq[T]$ of degree $k+d$ in terms of the polynomials $H_d$.

\begin{proposition}  \label{prop: formula para Hf}
Let $f:=T^{k+d}+f_{d-1}T^{k+d-1}+\cdots+f_0T^k$ be a polynomial of
$\fq[T]$ and let $H_f\in\fq[X_1,\ldots,X_{k+1}]$ be the polynomial
associated to $f$. Then the following identity holds:
\begin{equation} \label{eq:formula recursiva para H_d}
H_f=H_d+f_{d-1}H_{d-1}+\cdots+f_1H_1+f_0.
\end{equation}
\end{proposition}
\begin{proof}
In the proof of Lemma \ref{lema: formula recursiva para los H0d} we
obtain the following congruence relation:
$$T^{k+d}\equiv \Pi_1T^{k+d-1}-\Pi_2T^{k+d-2}+
\cdots+(-1)^{d-1}\Pi_dT^k+\mathcal{O}(T^{k-1})\mod{Q}.$$
Hence we have
$$
T^{k+d}+\sum_{j=0}^{d-1}f_jT^{k+j}\equiv
\sum_{j=0}^{d-1}\big((-1)^{d-1+j}\Pi_{d-j}+f_j\big)T^{k+j}+\mathcal{O}(T^{k-1})
\mod{Q}.$$
Therefore, taking into account that
$T^{k+j}\equiv H_jT^k+\mathcal{O}(T^{k-1})\mod{Q}$
holds for $1\leq j\leq d-1$, we obtain
$$
\begin{array}{rcl}
\displaystyle f:=T^{k+d}+\sum_{j=0}^{d-1}f_jT^{k+j}&\!\!\!\! \equiv
&\!\!\!\! \displaystyle \sum^{d-1}_{j=0}\big((-1)^{d-1+j}
\Pi_{d-j}+f_j\big)H_jT^k+\mathcal{O}(T^{k-1})\quad \mathrm{mod}\, Q\\[1ex]
  &\!\!\!\! = &\!\!\!\!\Bigg(\displaystyle
\displaystyle \sum^{d-1}_{j=0}(-1)^{d-1+j}\Pi_{d-j}H_j
+\sum^{d-1}_{j=0}f_jH_j\Bigg)T^k +
\mathcal{O}(T^{k-1})\\[1ex]
&\!\!\!\! = &\!\!\!\! \Bigg(H_d+\displaystyle
\sum^{d-1}_{j=0}f_jH_j\Bigg)T^k +\mathcal{O}(T^{k-1}),
\end{array}
$$
where the last equality is a consequence of Lemma \ref{lema: formula
recursiva para los H0d}. This shows that (\ref{eq:formula recursiva
para H_d}) is valid and finishes the proof.\end{proof}

\begin{remark}\label{observacion: acerca del grado de Hd entre otros}
From Lemma \ref{lema: formula recursiva para los H0d} and
Proposition \ref{Proposicion: formula explicita para la
hipersuperficie H_0d} we easily conclude that $H_d$ is a homogeneous
polynomial of $\fq[X_1,\ldots,X_{k+1}]$ of degree $d$ and can be
expressed as a polynomial in the elementary symmetric polynomials
$\Pi_1,\ldots,\Pi_d$. In this sense, we observe that $H_d$ is a
monic element of $\fq[\Pi_1,\ldots,\Pi_{d-1}][\Pi_d]$, up to a
nonzero constant of $\fq$. Combining these remarks and Proposition
\ref{prop: formula para Hf} we see that, for an arbitrary polynomial
$f:=T^{k+d}+f_{d-1}T^{k+d-1}+\cdots+f_0T^k\in\fq[T]$, the
corresponding polynomial $H_f\in \fq[X_1,\ldots,X_{k+1}]$ has degree
$d$ and is also a monic element of
$\fq[\Pi_1,\ldots,\Pi_{d-1}][\Pi_d]$.
 \end{remark}

\section{The geometry of the set of zeros of $H_f$}
\label{section:geometry of Vf}
For positive integers $d$ and $k$ with $k>d$, let be given
$f:=T^{k+d}+f_{d-1}T^{k+d-1}+\cdots+ f_0T^k\in\fq[T]$ and consider
the corresponding polynomial $H_f\!\in\!\fq[X_1,\ldots,X_{k+1}]$.
According to Remark \ref{observacion: acerca del grado de Hd entre
otros}, we may express $H_f$ as a polynomial in the first $d$
elementary symmetric polynomials $\Pi_1,\ldots,\Pi_d$ of
$\fq[X_1,\ldots,X_{k+1}]$, namely $H_f=G_f(\Pi_1,\ldots,\Pi_d)$,
where $G_f\in \fq[Y_1,\ldots,Y_d]$ is a monic element of
$\fq[Y_1,\ldots,Y_{d-1}][Y_d]$ of degree $1$ in $Y_d$.

In this section we obtain critical information on the geometry of
the set of zeros of $H_f$ that will allow us to establish upper
bounds on the number $q$--rational zeros of $H_f$.

\subsection{Notions of algebraic geometry}
Since our approach relies heavily on tools of algebraic geometry, we
briefly collect the basic definitions and facts that we need in the
sequel.  We use standard notions and notations of algebraic
geometry, which can be found in, e.g., \cite{Kunz85},
\cite{Shafarevich94}.

We denote by $\A^n$ the affine $n$--dimensional space $\cfq{\!}^{n}$
and by $\Pp^n$ the projective $n$--dimensional space over
$\cfq{\!}^{n+1}$. Both spaces are endowed with their respective
Zariski topologies, for which a closed set is the zero locus of
polynomials of $\cfq[X_1,\ldots, X_{n}]$ or of homogeneous
polynomials of  $\cfq[X_0,\ldots, X_{n}]$. For $\K:=\fq$ or
$\K:=\cfq$, we say that a subset $V\subset \A^n$ is an {\sf affine
$\K$--variety} if it is the set of common zeros in $\A^n$ of
polynomials $F_1,\ldots, F_{m} \in \K[X_1,\ldots, X_{n}]$.
Correspondingly, a {\sf projective $\K$--variety} is the set of
common zeros in $\Pp^n$ of a family of homogeneous polynomials
$F_1,\ldots, F_m \in\K[X_0 ,\ldots, X_n]$. An affine or projective
$\K$--variety is sometimes called simply a variety. When $V$ is the
set of zeros of a single polynomial of $\K[X_1,\ldots, X_{n}]$, or a
single homogeneous polynomial of $\K[X_0 ,\ldots, X_n]$, we say that
$V$ is an (affine or projective) {\sf $\fq$--hypersurface}.

A $\K$--variety $V$ is $\K$--{\sf irreducible} if it cannot be
expressed as a finite union of proper $\K$--subvarieties of $V$.
Further, $V$ is {\sf absolutely irreducible} if it is irreducible as
a $\cfq$--variety. An $\fq$--hypersurface $V$ is absolutely
irreducible if and only if any polynomial of $\fq[X_1,\ldots,
X_{n}]$, or any homogeneous polynomial of $\fq[X_0 ,\ldots, X_n]$,
of minimal degree defining $V$ is absolutely irreducible, namely is
an irreducible element of the ring $\cfq[X_1,\ldots,X_n]$ or
$\cfq[X_0 ,\ldots, X_n]$. Any $\K$--variety $V$ can be expressed as
an irredundant union $V=\mathcal{C}_1\cup \cdots\cup\mathcal{C}_s$
of absolutely irreducible $\K$--varieties, unique up to reordering,
which are called the {\sf absolutely irreducible} $\K$--{\sf
components} of $V$.

The set $V(\fq):=V\cap \fq^n$ is the set of {\sf $q$--rational
points} of $V$. Studying the number of elements of $V(\fq)$ is a
classical problem. The existence of $q$--rational points depends
upon many circumstances concerning the geometry of the underlying
variety.

For a $\K$-variety $V$ contained in $\A^n$ or $\Pp^n$, we denote by
$I(V)$ its {\sf defining ideal}, namely the set of polynomials of
$\K[X_1,\ldots, X_n]$, or of $\K[X_0,\ldots, X_n]$, vanishing on
$V$. The {\sf coordinate ring} $\K[V]$ of $V$ is the quotient ring
$\K[X_1,\ldots,X_n]/I(V)$ or $\K[X_0,\ldots,X_n]/I(V)$. The {\sf
dimension} $\dim V$ of a $\K$-variety $V$ is the length $r$ of the
longest chain $V_0\varsubsetneq V_1 \varsubsetneq\cdots
\varsubsetneq V_r$ of nonempty irreducible $\K$-varieties contained
in $V$. The {\sf degree} $\deg V$ of an irreducible $\K$-variety $V$
is the maximum number of points lying in the intersection of $V$
with a generic linear space $L$ of codimension $\dim V$, for which
$V\cap L$ is a finite set. More generally, following \cite{Heintz83}
(see also \cite{Fulton84}), if $V=V_1\cup\cdots\cup V_s$ is the
decomposition of $V$ into irreducible $\K$--components, we define
the degree of $V$ as
$$\deg V:=\sum_{i=1}^s\deg V_i.$$

Let $V$ be a variety contained in $\A^n$ and let $I(V)\subset
\cfq[X_1,\ldots, X_n]$ be the defining ideal of $V$. Let $\mathbf{x}$ be
a point of $V$. The {\sf dimension} $\dim_\mathbf{x}V$ {\sf of} $V$
{\sf at} $\mathbf{x}$ is the maximum of the dimensions of the
irreducible components of $V$ that contain $\mathbf{x}$. If
$I(V)=(F_1,\ldots, F_m)$, the {\sf tangent space}
$\mathcal{T}_\mathbf{x}V$ to $V$ at $\mathbf{x}$ is the kernel of
the Jacobian matrix $(\partial F_i/\partial X_j)_{1\le i\le m,1\le
j\le n}(\mathbf{x})$ of the polynomials $F_1,\ldots, F_m$ with respect
to $X_1,\ldots, X_n$ at $\mathbf{x}$. The point $\mathbf{x}$ is {\sf
regular} if $\dim\mathcal{T}_\mathbf{x}V=\dim_\mathbf{x}V$ holds.
Otherwise, the point $\mathbf{x}$ is called {\sf singular}. The set
of singular points of $V$ is the {\sf singular locus}
$\mathrm{Sing}(V)$ of $V$. For a projective variety, the concepts of
tangent space, regular and singular point can be defined by
considering an affine neighborhood of the point under consideration.

\subsection{The singular locus of symmetric hypersurfaces}
\label{subsec: singular locus symmetric hypersurfaces} With the
notations of the beginning of Section \ref{section:geometry of Vf},
let $V_f\subset \A^{k+1}$ denote the $\fq$--hypersurface defined by
$H_f$. Our main concern in this section is the study of the singular
locus of $V_f$. For this purpose, we consider the somewhat more
general framework that we now introduce. This will allow us to make
more transparent the facts concerning the algebraic structure of the
family of polynomials $H_f$ which are important at this point.

Let $Y_1,\ldots, Y_d$ be new indeterminates over $\cfq$, let
$G\in\fq[Y_1,\ldots,Y_d]$ be a given polynomial and let $\nabla G\in
\fq[Y_1,\ldots,Y_d]^d$ denote the vector consisting of the first
partial derivatives of $G$. Suppose that $\nabla G(\mathbf{y})$ is a
nonzero vector of $\A^d$ for every $\mathbf{y}\in\A^{d}$. Hence $G$
is square--free and defines a nonsingular hypersurface
$W\subset\A^d$.

Let $\Pi_1,\ldots,\Pi_d$ be the elementary symmetric polynomials of
$\fq[X_1,\ldots, X_{k+1}]$ and let $H:=G(\Pi_1,\ldots, \Pi_d)$. We
denote by $V\subset\A^{k+1}$ the hypersurface defined by $H$. The
main result of this section will be an upper bound on the dimension
of the singular locus of $V$. For this purpose, we consider the
following surjective morphism of $\fq$--hypersurfaces:
 \begin{eqnarray*}
   \Pi : V & \rightarrow & W
   \\
   \mathbf{x} & \mapsto & (\Pi_1(\mathbf{x}),\ldots,\Pi_d(\mathbf{x})).
 \end{eqnarray*}
For $\mathbf{x}\in V$ and $\mathbf{y}:=\Pi(\mathbf{x})$, we
denote by $\mathcal{T}_\mathbf{x} V$ and $\mathcal{T}_{\mathbf{y}}
W$ the tangent spaces to $V$ at $\mathbf{x}$ and to $W$ at
$\mathbf{y}$. We also consider the differential map of $\Pi$
at $\mathbf{x}$, namely
 \begin{eqnarray*}
   \mathrm{d}_\mathbf{x}\Pi :  \mathcal{T}_\mathbf{x} V &
   \rightarrow & \mathcal{T}_{\mathbf{y}}W \\
   \mathbf{v} & \mapsto & A(\mathbf{x})\cdot \mathbf{v},
 \end{eqnarray*}
where $A(\mathbf{x})$ stands for the  $d \times (k+1)$ matrix
  \begin{equation} \label{eq:jacobiano de los pi}
    A(\mathbf{x}):=
       \left(
         \begin{array}{ccc}
           \dfrac{\partial\Pi_1}{\partial X_1}(\mathbf{x}) & \cdots & \dfrac{\partial\Pi_1}{\partial X_{k+1}}(\mathbf{x})
         \\
           \vdots & & \vdots
          \\
           \dfrac{\partial\Pi_d}{\partial X_1}(\mathbf{x})& \cdots & \dfrac{\partial\Pi_d}{\partial X_{k+1}}(\mathbf{x})
         \end{array}
       \right).
  \end{equation}
In order to prove our result about the singular locus of $V$, we
first make a few remarks concerning the Jacobian matrix of the
elementary  symmetric polynomials that will be useful in the sequel.

It is well known that the first partial derivatives of the
elementary symmetric polynomials $\Pi_i$  satisfy the following
equalities (see, e.g., \cite{LaPr02}) for $1\leq i,j \leq k+1$:
  \begin{equation}\label{eq:derivadas parciales simetricos elementales}
    \frac{\partial \Pi_i}{\partial X_{j}}= \Pi_{i-1}-X_{j} \Pi_{i-2}
    + X_{j}^2 \Pi_{i-3} +\cdots+ (-1)^{i-1} X_{j}^{i-1}.
  \end{equation}
As a consequence, denoting by $A_{k+1}$ the $(k+1)\times (k+1)$ Vandermonde matrix
  \begin{equation}\label{eq:matriz de vandermonde k+1}
    A_{k+1}:=\left(
         \begin{array}{cccc}
           1 & 1 & \cdots & 1
         \\
           X_1 & X_2 & \cdots &  X_{k+1}
          \\  \vdots & \vdots & & \vdots
          \\
           X_1^k & X_2^k & \cdots & X_{k+1}^k,
         \end{array}
       \right),
 \end{equation}
we deduce that the  Jacobian matrix of $\Pi_1,\dots,\Pi_{k+1}$ with
respect to $X_1,\dots,X_{k+1}$ can be factored as follows:
\begin{equation} \label{eq:factorizacion del jacobiano de los simetricos elementales}
     \left(\frac{\partial \Pi_i}{\partial X_j}\right)_{1\le i,j\le k+1}:=B_{k+1}\cdot
      A_{k+1}
       :=
        \left(
         \begin{array}{ccccc}
           1 & \ 0 & 0 &  \dots & 0
         \\
           \Pi_1 & - 1 & 0 &  &
          \\
           \Pi_2 & -\Pi_1 & 1 & \ddots & \vdots
          \\
           \vdots &\vdots  & \vdots & \ddots & 0
           \\
           \Pi_k & -\Pi_{k-1} & \Pi_{k-2} & \cdots & (-1)^{k}
         \end{array}
       \right)
     \cdot
         A_{k+1}
  \end{equation}
We observe that the left factor $B_{k+1}$ is a square,
lower--triangular matrix whose determinant is equal to
$(-1)^{{k(k+1)}/{2}}$. This implies that the determinant of the
Jacobian matrix $(\partial{\Pi_i}/\partial{X_j})_{1\leq i,j\leq
k+1}$ is equal, up to a sign, to the determinant of $A_{k+1}$, i.e.,
$$\det \left(\frac{\partial \Pi_i}{\partial X_j}\right)_{1\le i,j\le
k+1}=(-1)^{{k(k+1)}/{2}} \prod_{1\le i < j\le k+1}(X_i-X_j).$$
An interesting fact, which will not be used in what follows, is that
the inverse matrix of the matrix $B_{k+1}$ of (\ref{eq:factorizacion
del jacobiano de los simetricos elementales}) is given by
%
$$B_{k+1}^{-1}=\left(
         \begin{array}{ccccc}
           H_0 & 0 & 0 &  \dots & 0
         \\
           H_1 & - H_0 & 0 &   &
          \\
           H_2 & -H_1 & H_0 &  \ddots & \vdots
          \\
           \vdots & \vdots &\vdots & \ddots & 0
           \\
           H_k & -H_{k-1} & H_{k-2} & \cdots & (-1)^{k}H_0
         \end{array}
       \right).$$

\begin{theorem}\label{teo:caracterizacion de los puntos singulares}
  The singular locus $\Sigma$ of $V$ has  dimension  at most $d-1$.
  Moreover, the elements of $\Sigma$ have at most $d-1$
  pairwise--distinct coordinates.
\end{theorem}
\begin{proof}
    By the chain rule we deduce that the partial derivatives of $H$
    satisfy the following equality for $1\leq j \leq k+1$:
    $$
      \dfrac{\partial H}{\partial X_j}  =    \bigg(\dfrac{\partial G}{\partial Y_1} \circ \Pi \bigg)
      \cdot\dfrac{\partial\Pi_1}{\partial X_j} + \cdots +
      \bigg(\dfrac{\partial G}{\partial Y_d}\circ \Pi \bigg)
      \cdot \dfrac{\partial\Pi_d}{\partial X_j}
    $$
    If $\mathbf{x}$ is any  point of $\Sigma$, then we have
    $$
    \nabla H (\mathbf{x})=\nabla G (\Pi (\mathbf{x}))\cdot A (\mathbf{x})= \mathbf{0},
    $$
where $A(\mathbf{x})$ is the matrix  defined in (\ref{eq:jacobiano
de los pi}). Fix $\mathbf{x}\in \Sigma$ and let $\mathbf{v}: =
\nabla G\big(\Pi(\mathbf{x})\big)$. By hypothesis we have that
$\mathbf{v} \in \A^{d}$ is a nonzero vector satisfying
    $$\mathbf{v}\cdot A(\mathbf{x})= \mathbf{0}.$$
Hence, all the maximal minors of $A(\mathbf{x})$ must be zero.

The matrix $A(\mathbf{x})$ is the $d\times (k+1)$--submatrix of
$(\partial{\Pi_i}/\partial{X_j})_{1\leq i,j\leq k+1}(\mathbf{x})$
consisting of the first $d$ rows of the latter. Therefore, from
(\ref{eq:factorizacion del jacobiano de los simetricos elementales})
we conclude that
$$
A(\mathbf{x})=B_{d,k+1}(\mathbf{x})\cdot A_{k+1}(\mathbf{x}),
$$
where $B_{d,k+1}(\mathbf{x})$ is the $d\times (k+1)$--submatrix of
$B_{k+1}(\mathbf{x})$ consisting of the first $d$ rows of
$B_{k+1}(\mathbf{x})$. Furthermore, since the last $k+1-d$ columns
of $B_{d,k+1}(\mathbf{x})$ are zero, we may rewrite this identity in
the following way:
\begin{equation}\label{eq:factorization submatrix Jacobian elem symm}
A(\mathbf{x})=B_d(\mathbf{x})\cdot \left(
         \begin{array}{cccc}
           1 & 1 & \dots & 1
         \\
           x_1 & x_2 & \dots &  x_{k+1}
          \\  \vdots & \vdots & & \vdots
          \\
           x_1^{d-1} & x_2^{d-1} & \dots & x_{k+1}^{d-1},
         \end{array}
       \right),
\end{equation}
where $B_d(\mathbf{x})$ is the $(d\times d)$--submatrix of
$B_{k+1}(\mathbf{x})$ consisting on the first $d$ rows and the first
$d$ columns of $B_{k+1}(\mathbf{x})$.

Fix $1\leq l_1<\cdots<l_d\leq k+1$, set $I:=(l_1,\ldots,l_d)$ and
consider the $(d\times d)$--submatrix $M_I(\mathbf{x})$ of
$A(\mathbf{x})$ consisting of the columns $l_1,\ldots,l_d$ of
$A(\mathbf{x})$, namely $M_I(\mathbf{x}):=({\partial \Pi_i}/\partial
X_{l_{j}})_{1\leq i,j\leq d}(\mathbf{x})$.

From (\ref{eq:factorizacion del jacobiano de los simetricos
elementales}) and (\ref{eq:factorization submatrix Jacobian elem
symm}) we easily see that $M_I(\mathbf{x})=B_d(\mathbf{x})\cdot
A_{d,I}(\mathbf{x})$, where $A_{d,I}(\mathbf{x})$ is the Vandermonde
matrix $A_{d,I}(\mathbf{x}):=(x_{l_{j}}^{i-1})_{1\leq i,j\leq d}$.
Therefore, we obtain
    \begin{equation}\label{eq:determinante de vandermonde}
      \det\big(M_I(\mathbf{x})\big) = (-1)^{{(d-1)d}/{2}}\det A_{d,I}(\mathbf{x})
                = (-1)^{{(d-1)d}/{2}} \prod_{1\le r<s\le d}(x_{l_{r}}-x_{l_{s}})=0.
    \end{equation}

    Since (\ref{eq:determinante de vandermonde}) holds for every
    $I:=(l_1,\ldots,l_d)$ as above, we conclude that every point
    $\mathbf{x}\in \Sigma$ has at most $d-1$ pairwise--distinct
    coordinates. In particular, $\Sigma$ is contained in a finite union
    of linear varieties of $\A^{k+1}$ of dimension $d-1$, and thus its
    dimension is at most $d-1$.
\end{proof}

We observe that the proof of Theorem \ref{teo:caracterizacion de los
puntos singulares} provides a more precise description of the
singular locus $\Sigma$ of $V$, which is the subject of the
following remark.

 \begin{remark}\label{obs: caracterizacion de las variedades lineales}
  Let notations and assumptions be as in Theorem
  \ref{teo:caracterizacion de los puntos singulares}. From the proof
  of Theorem \ref{teo:caracterizacion de los puntos singulares} we
  obtain the following inclusion:
  $$
    \Sigma\subset \bigcup_{\mathcal{I}}\mathcal{L_{I}},
  $$
where $\mathcal{I}:=\{I_1,\ldots,I_{d-1}\}$ runs over all the
partitions of $\{1,\ldots,k+1\}$ into $d-1$ nonempty subsets
$I_j\subset \{1,\ldots,k+1\}$ and $\mathcal{L_{I}}$ is the linear
variety
  $$
    \mathcal{L_{I}}:= \mathrm{span}(\mathbf{v}^{(I_1)},\ldots,\mathbf{v}^{(I_{d-1})})
  $$
spanned by the vectors $\mathbf{v}^{(I_j)}:=
(v_1^{(I_j)},\ldots,v_{k+1}^{(I_j)})$ defined by $v_m^{(I_j)}:=1$
for $m\in I_j$ and $v_m^{(I_j)}:=0$ for $m\notin I_j$. In
particular, it follows that if $\Sigma$ has dimension $d-1$, then it
contains a linear variety $\mathcal{L_{I}}$ as above.
\end{remark}

\subsection{The dimension of the singular locus of $V_f$}

Now we consider the hypersurface $V_f$ defined by the polynomial
$H_f\in \fq[X_1,\ldots X_{k+1}]$ associated to the polynomial
$f:=T^{k+d}+f_{d-1}T^{k+d-1}+\cdots+ f_0T^k$. According to Remark
\ref{observacion: acerca del grado de Hd entre otros}, we may
express $H_f$ in the form $H_f=G_f(\Pi_1,\ldots,\Pi_d)$, where
$G_f\in \fq[Y_1,\ldots,Y_d]$ is a polynomial of degree $d$ which is
monic in $Y_d$, up to a nonzero constant. Moreover, since
$$
  \nabla G_f (\mathbf{y}) = \left(\dfrac{\partial G_f}{\partial Y_1}(\mathbf{y})
  ,\ldots, \dfrac{\partial G_f}{\partial Y_d}(\mathbf{y}), (-1)^{d-1}\right)
$$
holds for every $\mathbf{y}\in \mathbb{A}^d$, we see that $\nabla
G_f(\mathbf{y})\not=\mathbf{0}$ for every $\mathbf{y}\in\A^d$; in
other words, $G_f$ defines a nonsingular hypersurface
$W\subset\A^d$. Then the results of Section 3.2  can be applied to
$H_f$. In particular, we have the following immediate consequence of
Theorem \ref{teo:caracterizacion de los puntos singulares}.
\begin{corollary}\label{coro:dimension singular locus Vf}
The singular locus $\Sigma_f\subset\A^{k+1}$ of $V_f$ has dimension
at most $d-1$.
\end{corollary}

In order to obtain estimates on the number of $q$--rational points
of $V_f$ we also need information concerning the behavior of $V_f$
``at infinity''. For this purpose, we consider the projective
closure $\mathrm{pcl}(V_f)\subset\mathbb{P}^{k+1}$ of $V_f$, whose
definition we now recall. Consider the embedding of $\A^{k+1}$ into
the projective space $\Pp^{k+1}$ which assigns to any
$\mathbf{x}:=(x_1,\ldots, x_{k+1})\in\A^{k+1}$ the point
$(1:x_1:\dots:x_{k+1})\in\Pp^{k+1}$. The closure
$\mathrm{pcl}(V_f)\subset\Pp^{k+1}$ of the image of $V_f$ under this
embedding in the Zariski topology of $\Pp^{k+1}$ is called the {\sf
projective closure} of $V_f$. The points of $\mathrm{pcl}(V_f)$
lying in the hyperplane $\{X_0=0\}$ are called the points of
$\mathrm{pcl}(V_f)$ at infinity.

It is well--known that $\mathrm{pcl} (V_f)$ is the
$\fq$--hypersurface of $\mathbb{P}^{k+1}$ defined by the
homogenization $H_f^h\in\fq[X_0,\ldots,X_{k+1}]$ of the polynomial
$H_f$ (see, e.g., \cite[\S I.5, Exercise 6]{Kunz85}). We have the
following result.
\begin{proposition}\label{prop:dimension singularidades en el infinito}
$\mathrm{pcl}(V_f)$ has singular locus at infinity of dimension  at
most $d-2$.
\end{proposition}
\begin{proof}
By Proposition \ref{prop: formula para Hf}, we have
$$H_f=H_d+f_{d-1}H_{d-1}+\cdots+f_1H_1+f_0,$$
where each $H_j$ is a homogeneous polynomial of degree $j$ for
$1\leq j\leq d$. Hence, the homogenization of $H_f$ is the following
polynomial of $\fq[X_0,\ldots, X_{k+1}]$:
\begin{equation} \label{eq:homogenizacion de H_d}
H_f^h=H_d+f_{d-1}H_{d-1}X_0+\cdots+f_1H_1X_0^{d-1}+f_0X_0^d.
\end{equation}
Let $\Sigma_f^{\infty}\subset\mathbb{P}^{k+1}$ denote the singular
locus of $\mathrm{pcl}(V_f)$ at infinity, namely the set of singular
points of $\mathrm{pcl}(V_f)$ lying in the hyperplane $\{X_0=0\}$.
We have that any point $\mathbf{x}\in\Sigma_f^{\infty}$ satisfies
the identities $H_f^h(\mathbf{x})=0$ and ${\partial H_f^h}/{\partial
X_i}(\mathbf{x})=0$ for $0\leq i \leq k+1$. From
(\ref{eq:homogenizacion de H_d}) we see that any point
$\mathbf{x}:=(0:x_1:\cdots:x_{k+1})\in \Sigma_f^{\infty} $ satisfies
the identities
\begin{equation}\label{sistema de las derivadas en la homogenizacion}
  \left\{\begin{array}{rlll}
     H_d(x_1,\ldots,x_{k+1})&=& 0 ,&  \\[1ex]
      f_{d-1}H_{d-1} (x_1,\ldots,x_{k+1})&=& 0 ,&  \\[1ex]
\displaystyle
    \frac{\partial H_d}{\partial X_i}(x_1,\ldots,x_{k+1})&=& 0 & (1\leq i\leq
    k+1).
   \end{array}\right.
\end{equation}
From Proposition \ref{Proposicion: formula explicita para la
hipersuperficie H_0d} and Remark \ref{observacion: acerca del grado
de Hd entre otros} we have that $H_d\in \fq[X_1,\ldots,X_{k+1}]$ is
a homogeneous polynomial of degree $d$ which can be expressed in the
form $H_d=G_d(\Pi_1,\ldots,\Pi_d)$, where
$G_d\in\fq[Y_1,\ldots,Y_d]$ has degree $d$ and is monic in $Y_d$.
Combining these remarks with Theorem \ref{teo:caracterizacion de los
puntos singulares} we conclude that the set of solutions of
(\ref{sistema de las derivadas en la homogenizacion}) is an affine
cone of $\A^{k+1}$ of dimension at most $d-1$, and hence, a
projective variety of $\Pp^{k}$ of dimension at most $d-2$. This
finishes the proof of the proposition.
\end{proof}

We end this section with a useful consequence of our bound on the
dimension of the singular locus of $\mathrm{pcl}(V_f)$, namely that
$V_f$ is absolutely irreducible. This result, which has been proved
in \cite[Section 4]{ChMu07}, is obtained here as an easy consequence
of Proposition \ref{prop:dimension singularidades en el infinito}.
\begin{corollary}\label{coro:Vf is abs irred} The hypersurface $V_f$ is
absolutely irreducible.
\end{corollary}
\begin{proof}
We observe that $V_f$ is absolutely irreducible if and only if
$\mathrm{pcl}(V_f)$ is absolutely irreducible (see, e.g.,
\cite[Chapter I, Proposition 5.17]{Kunz85}). If $\mathrm{pcl}(V_f)$
is not absolutely irreducible, then it has a nontrivial
decomposition into absolutely irreducible components
$$\mathrm{pcl}(V_f)=\mathcal{C}_1\cup\cdots\cup \mathcal{C}_s,$$
where $\mathcal{C}_1,\ldots,\mathcal{C}_s$ are projective hypersurfaces
of $\Pp^{k+1}$. Since $\mathcal{C}_i\cap\mathcal{C}_j\neq\emptyset$
and $\mathcal{C}_i$, $\mathcal{C}_j$ are absolutely irreducible, we
conclude that $\dim (\mathcal{C}_i\cap \mathcal{C}_j)=k-1$ holds.

Denote by $\Sigma_f^{h}$ the singular locus of $\mathrm{pcl}(V_f)$.
Corollary \ref{coro:dimension singular locus Vf} and Proposition
\ref{prop:dimension singularidades en el infinito} imply
$\dim\Sigma_f^h\le d-1$. On the other hand, we have
$\mathcal{C}_i\cap \mathcal{C}_j\subset\Sigma_f^h$ for any $i\not=
j$, which implies $\dim \Sigma_f^h\geq k-1$. This contradicts the
assertion $\dim\Sigma_f^{h}\leq d-1$, since we have $d<k$ by
hypothesis. It follows that $V_f$ is absolutely irreducible.
\end{proof}

\section{The singular locus of $V_f$ for fields of large characteristic}
\label{section:singular locus for large characteristic}
In this section we characterize the set of polynomials $f\in\fq[T]$
for which the associated hypersurface $V_f\subset \mathbb{A}^{k+1}$
has a singular locus of dimension $d-1$. This characterization
enables us giving conditions under which such polynomials do not
generate deep holes of the standard Reed--Solomon code of dimension
$k$ over $\fq$.

The first step is to obtain a suitable expression of the derivatives
of the polynomial $H_d$ associated to $T^{k+d}$.
\begin{lemma} \label{Lema: derivadas de los H0d}
If $j\geq 2$, then the partial derivatives of the polynomials $H_j$
satisfy the following identity for $1\leq i \leq k+1$:
 $$ 
 \frac{\partial H_j}{\partial X_i}=H_{j-1}+H_{j-2}X_i+\cdots+H_{j-3}X_i^2+\cdots + X_i^{j-1}.
 $$
\end{lemma}
\begin{proof}
The proof is by induction on $j$. By Lemma \ref{lema: formula
recursiva para los H0d} we have $H_1=\Pi_1$ and $H_2= \Pi_1H_1 -
\Pi_2$. Combining this with (\ref{eq:derivadas parciales simetricos
elementales}) we easily see the assertion for $j=2$. Next assume
that the statement of the lemma holds for $2\leq j \leq l-1$; we
claim that it also holds for $j=l$. According to Lemma \ref{lema:
formula recursiva para los H0d}, we have
\begin{equation} \label{eq:Lema derivadas de los H0d, primera sumatoria}
\dfrac{\partial H_l}{\partial X_i}=\sum_{m=1}^l
(-1)^{m-1}\dfrac{\partial (\Pi_m H_{l-m})}{\partial X_i}.
\end{equation}
By the inductive hypothesis and the expression (\ref{eq:derivadas
parciales simetricos elementales}) for the first partial derivatives
of the elementary symmetric polynomials, each term in the
right--hand side of (\ref{eq:Lema derivadas de los H0d, primera
sumatoria}) can be expressed as follows:
\begin{equation}
\label{eq:lema derivadas de los H0d segunda formula de la
demostracion} \dfrac{\partial(\Pi_m H_{l-m})}{\partial X_i}=
H_{l-m}\sum_{n=1}^m(-1)^{n-1}\Pi_{m-n}X_i^{n-1}+
\Pi_m\sum_{n=1}^{l-m}H_{l-(m+n)}X_i^{n-1}.
\end{equation}
Now we determine the coefficient of $H_{l-m}$ in the right--hand
side of (\ref{eq:Lema derivadas de los H0d, primera sumatoria}).
From (\ref{eq:lema derivadas de los H0d segunda formula de la
demostracion}) we see that the only terms having a nonzero
contribution to the coefficient of $H_{l-m}$ are ${\partial (\Pi_n
H_{l-n})}/{\partial X_i}$ for $1\leq n \leq m$. In particular, we
easily deduce that the coefficient of $H_{l-1}$ is $1$. For $1\le n
<m$, the summand $(-1)^{n-1}{\partial(\Pi_n H_{l-n})}/{\partial
X_i}$ contributes with the term $(-1)^{n-1}X_i^{m-n-1}\Pi_n$. On the
other hand, the summand $(-1)^{m-1}{\partial(\Pi_m
H_{l-m})}/{\partial X_i}$ in the left--hand side of (\ref{eq:lema
derivadas de los H0d segunda formula de la demostracion})
contributes with the sum
$(-1)^{m-1}\sum_{n=0}^{m-1}(-1)^{m-n-1}\Pi_n X_i^{m-n-1}.$ Putting
all these terms together, we conclude that the term $H_{l-m}$ occurs
in (\ref{eq:Lema derivadas de los H0d, primera sumatoria})
multiplied by
$$(-1)^{m-1}\sum_{n=0}^{m-1}(-1)^{m-n-1}\Pi_n X_i^{m-n-1}+
\sum_{n=1}^{m-1}(-1)^{n-1}\Pi_n X_i^{m-n-1}=X_i^{m-1}.$$
This finishes the proof of the lemma.
\end{proof}

Observe that, similarly to the factorization (\ref{eq:factorizacion
del jacobiano de los simetricos elementales}) of the Jacobian matrix
of the elementary symmetric polynomials of Section \ref{subsec:
singular locus symmetric hypersurfaces}, Lemma \ref{Lema: derivadas
de los H0d} allows us to express the Jacobian matrix of $H_1,\dots,
H_{k+1}$ with respect to $X_1,\dots,X_{k+1}$ as the following matrix
product:
 \begin{equation} \label{eq:factorizacion del jacobiano de los Hi}
     \left(\frac{\partial H_i}{\partial X_j}\right)_{1\leq i, j \leq k+1}
     :=
      \left(
         \begin{array}{ccccc}
           H_0 & 0 &   \cdots & 0
         \\
           H_1 &  H_0 & \ddots   & \vdots
          \\
           \vdots & \vdots & \ddots &0
           \\
           H_k & H_{k-1} & \dots & H_0
         \end{array}
       \right)
     \cdot
         A_{k+1},
  \end{equation}
where $A_{k+1}$ is the Vandermonde matrix defined in (\ref{eq:matriz de vandermonde k+1}).

Let $f\in\fq[T]$ be a polynomial of the form
 $$
   f:=T^{k+d}+f_{d-1}T^{k+d-1}+\cdots+ f_1T^{k+1}+f_0T^k,
 $$
and let $V_f\subset \mathbb{A}^{k+1}$ be the hypersurface associated
to $f$. By Proposition \ref{prop: formula para Hf}, we have that
$V_f$ is the hypersurface defined by the polynomial
 $$
   H_f= H_d+f_{d-1}H_{d-1}+\cdots+ f_1H_1+f_0 H_0,
 $$
where the polynomials $H_j\in\fq[X_1,\ldots,X_{k+1}]$ $(0\leq j \leq
d)$ are defined in Section \ref{section:planteo del problema}. We
recall that each $H_j$ is homogeneous and symmetric of degree $j$
(Remark \ref{observacion: acerca del grado de Hd entre otros}).

Corollary  \ref{coro:dimension singular locus Vf} asserts that the
singular locus $\Sigma_f$ of $V_f$ has dimension at most $d-1$.
Suppose now that $\Sigma_f$ has dimension $d-1$. From Remark
\ref{obs: caracterizacion de las variedades lineales} we see that
there exists a partition $\mathcal{I}:=\{I_1,\ldots,I_{d-1}\}$ of
the set $\{1,\ldots,k+1\}$ into $d-1$ nonempty sets
$I_j\subset\{1,\ldots,k+1\}$ with the following property: let
$\mathcal{L_I}\subset \A^{k+1}$ denote the linear variety
  $$
    \mathcal{L_I}:=\mathrm{span}(\mathbf{v}^{(I_1)},\ldots, \mathbf{v}^{(I_{d-1})})
  $$
spanned by the vectors
$\mathbf{v}^{(I_j)}:=(v_1^{(I_j)},\ldots,v_{k+1}^{(I_{j})})$ defined
by $v_l^{(I_j)}:=1$ for $l\in I_j$ and $v_l^{(I_j)}:=0$ for $l\notin
I_j$ $(1\leq j \leq d-1)$. Then  $\mathcal{L}_I\subset\Sigma_f$
holds. Let $\lambda:=(\lambda_1,\ldots,\lambda_{d-1})\in \A^{d-1}$
and let $\mathbf{x}:=\sum_{j=1}^{d-1}\lambda_j \mathbf{v}^{(I_j)}$
be an arbitrary point of $\mathcal{L_I}$. Since $\mathbf{x}$ is a
singular point of $V_f$ we have
  $$
    0=\dfrac{\partial H_f}{\partial X_i}(\mathbf{x})= \dfrac{\partial H_d}{\partial X_i}
    (\mathbf{x})+ \sum_{j=1}^{d-1}f_{d-j}\dfrac{\partial H_{d-j}}{\partial X_i}(\mathbf{x})
  $$
for $1\leq i \leq k+1$. This shows that the following matrix identity holds:
\begin{equation}\label{eq:singular locus de dimension d-1 primer sistema }
  -\left(%
    \begin{array}{ccc}
      \dfrac{\partial H_1}{\partial X_1}(\mathbf{x}) & \cdots
      & \dfrac{\partial H_{d-1}}{\partial X_1}(\mathbf{x})
      \\
       \vdots &  & \vdots
      \\
       \dfrac{\partial H_1}{\partial X_{k+1}}(\mathbf{x}) & \cdots
       & \dfrac{\partial H_{d-1}}{\partial X_{k+1}}(\mathbf{x})
      \\
    \end{array}%
   \right)
   \left(%
    \begin{array}{c}
      f_1  \\  \vdots   \\ f_{d-1}
    \end{array}%
   \right)
  =
   \left(%
    \begin{array}{c}
       \dfrac{\partial H_d}{\partial X_1}(\mathbf{x})
      \\
       \vdots
      \\
       \dfrac{\partial H_d}{\partial X_{k+1}}(\mathbf{x}) \\
    \end{array}%
  \right).
\end{equation}

By symmetry, we may assume that $x_i=\lambda_i$ holds for $1\leq i
\leq d-1$. We further assume that $\lambda_i\neq \lambda_j$ for
$i\neq j$. Considering the first $d-1$ equations of
(\ref{eq:singular locus de dimension d-1 primer sistema }) we obtain
the square system
  \begin{equation}\label{eq:singular locus de dimension d-1 primeras  d-1 ecuaciones}
    -B(\mathbf{x})
     \cdot \left(%
        \begin{array}{c}
          f_1  \\  \vdots  \\ f_{d-1}
        \end{array}%
       \right)
      =
       \left(%
        \begin{array}{c}
           \dfrac{\partial H_d}{\partial X_1}(\mathbf{x})
          \\
           \vdots
          \\
           \dfrac{\partial H_d}{\partial X_{d-1}}(\mathbf{x}) \\
        \end{array}%
       \right),
  \end{equation}
where $B(\mathbf{x})\in\A^{(d-1)\times (d-1)}$ is the matrix
  $$
    B(\mathbf{x}):=\left(%
     \begin{array}{ccc}
      \dfrac{\partial H_1}{\partial X_1}(\mathbf{x}) & \cdots
      & \dfrac{\partial H_{d-1}}{\partial X_1}(\mathbf{x})
      \\
       \vdots &  & \vdots
      \\
       \dfrac{\partial H_1}{\partial X_{d-1}}(\mathbf{x}) & \cdots
       & \dfrac{\partial H_{d-1}}{\partial X_{d-1}}(\mathbf{x})
      \\
    \end{array}%
   \right).
  $$
From (\ref{eq:factorizacion del jacobiano de los Hi}) we see that
$B(\mathbf{x})$ can be factored as follows:
    \begin{equation}\label{eq:factorizacion de B(x)}
      B(\mathbf{x})
   =    A_{d-1}(\mathbf{x})^t
       \cdot
   \left(
         \begin{array}{ccccc}
           H_0 & H_1 (\mathbf{x})\!\!\!\! &   \dots & H_{d-2}(\mathbf{x})
         \\
           0 &  H_0 &  \dots & H_{d-3}(\mathbf{x})
          \\
           \vdots & \ddots  & \ddots & \vdots
           \\
           0 &  \dots & 0 & H_0
         \end{array}
       \right),
  \end{equation}
where $A_{d-1}(\mathbf{x})$ is the Vandermonde matrix
$A_{d-1}(\mathbf{x}):=(x_j^{i-1})_{1\leq i, j \leq d-1}$. As a
consequence, we have that $B(\mathbf{x})$ is nonsingular and its
determinant is equal to
  $$
    \det B(\mathbf{x}) = \prod_{1\leq i<j\leq d-1}(x_i-x_j).
  $$
Hence $(f_1,\ldots,f_{d-1})$ is the unique solution of the linear
system (\ref{eq:singular locus de dimension d-1 primeras  d-1
ecuaciones}). Furthermore, by the Cramer rule we obtain
  $$
    f_j=\dfrac{\det B^{(j)}(\mathbf{x})}{\det B(\mathbf{x})}\quad(1\leq j \leq d-1),
  $$
where $B^{(j)}(\mathbf{x})\in\A^{(d-1)\times (d-1)}$ is the matrix
obtained by replacing the $j$th column of $B(\mathbf{x})$ by the
vector $b(\mathbf{x}):=\big((\partial H_d/\partial
X_1)(\mathbf{x}),\ldots,({\partial H_d}/{\partial
X_{d-1}})(\mathbf{x})\big)$.

Let $B,B^{(j)}\in\fq[X_1,\dots,X_{k+1}]^{(d-1)\times(d-1)}$ be the
``generic'' versions of the matrices
$B(\mathbf{x}),B^{(j)}(\mathbf{x})$ for $1\leq j \leq d-1$. We claim
that $\det B = \prod_{1\leq i<j\leq d-1}(X_i-X_j)$ divides $\det
B^{(j)}$ in $\fq[X_1, \ldots, X_{k+1}]$.

In order to show this claim, let $C\in\fq[X_1, \ldots,
X_{d-1}]^{(d-1)\times d}$ be the following matrix:
\begin{equation}\label{eq:matrix C}
   C:= \left(
    \begin{array}{ccccc}
     1 & X_1 & \cdots & X_{1}^{d-2} & X_{1}^{d-1}
     \\[1ex]
     1 & X_2 & \cdots & X_{2}^{d-2} & X_{2}^{d-1}
     \\ \vdots & \vdots &  & \vdots & \vdots
     \\
     1 & X_{d-1} & \cdots & X_{d-1}^{d-2} & X_{d-1}^{d-1}
    \end{array}
  \right).
\end{equation}
Observe that this matrix is obtained by appending the vector column
$(X_j^{d-1})_{1\le j\le d-1}$ to the generic matrix
$A_{d-1}\in\fq[X_1,\dots,X_{d-1}]^{(d-1)\times (d-1)}$. Further, for
$1\leq j \leq d-1$ we define a matrix $H^{(j)}\in
\fq[X_1,\dots,X_{k+1}]^{d\times (d-1)}$ as
 $$
   H^{(j)}
   :=\left(
    \begin{array}{cccccccc}
     H_0 & H_1 & \cdots & H_{j-2} & H_{d-1} & H_{j} & \cdots & H_{d-2}
     \\
      0 & H_0 & \cdots & H_{j-1} & H_{d-2}& H_{j-1} & \cdots &  H_{d-3}
     \\
      \vdots & 0 &  \ddots & \vdots & \vdots&\vdots&& \vdots
     \\
       &  & \ddots  &  H_0 & H_{d-j+1} & H_2 & \cdots & H_{j-1}
     \\
      \vdots & \vdots &   & 0 & H_{d-j} & H_1 &  & \vdots
      \\
       &  &  &  &H_{d-j-1}  & H_0 & \ddots &
      \\
       \vdots  & \vdots &  & \vdots &\vdots& 0&\ddots&H_1
      \\
       &  &  &  &H_1 & \vdots & \ddots & H_0
     \\
     0 & 0 & \cdots & 0 & H_0 & 0 & \cdots & 0
    \end{array}
  \right),
 $$
namely $H^{(j)}$ is obtained by appending a zero $d$th row to the
right factor in the right--hand side of (\ref{eq:factorizacion de
B(x)}) and replacing the resulting $j$th column by the column vector
$(H_{d-j}:1\le j\le d)\in\fq[X_1,\dots,X_{k+1}]^{d\times 1}$.

It turns out that the matrix $B^{(j)}$ can be factored as follows:
 \begin{equation}\label{eq:factorizacion de Bj}
    B^{(j)}= C  \cdot H^{(j)}
 \end{equation}
Indeed, for $l\not=j$, the $l$th columns of $B$ and $B^{(j)}$ agree,
and the fact that the $l$th columns of both sides of
(\ref{eq:factorizacion de Bj}) are equal is easily deduced from
(\ref{eq:factorizacion de B(x)}). On the other hand, from Lemma
\ref{Lema: derivadas de los H0d} we immediately conclude that the
$j$th columns of $B^{(j)}$ and $C\cdot H^{(j)}$ are equal.

In particular, the determinant of $B^{(j)}$ can be obtained from
(\ref{eq:factorizacion de Bj}) by means of the Cauchy--Binet
formula. Since any maximal minor of $C$ is a multiple of $\det B$
(see, e.g., \cite[Lemma 2.1]{Ernst00} or \cite[Exercise
281]{FaSo65}), we immediately deduce that $\det B$ divides $\det
B^{(j)}$ in $\fq[X_1, \ldots, X_{k+1}]$.

As a consequence of our claim, we see that for $1\leq j \leq d-1$
there exists a homogeneous polynomial
$P^{(j)}\in\fq[X_1,\ldots,X_{d-1}]$ of degree $d-j$ or zero such
that
 \begin{equation}\label{eq:seccion 4 f_d-j=P(j)}
  f_{d-j}=P^{(j)}(\lambda_1,\ldots,\lambda_{d-1})
 \end{equation}
holds for $1\leq j \leq d-1$ and for any $(\lambda_1,\ldots,
\lambda_{d-1})\in \A^{d-1}$ with $\lambda_i\neq \lambda_j$ for
$i\neq j$. Since (\ref{eq:seccion 4 f_d-j=P(j)}) holds in a Zariski
open dense subset of $\A^{d-1}$, we conclude that (\ref{eq:seccion 4
f_d-j=P(j)}) holds for every
$(\lambda_1,\ldots,\lambda_{d-1})\in\A^{d-1}$. By substituting $0$
for $\lambda_i$ in (\ref{eq:seccion 4 f_d-j=P(j)}) we deduce that
$f_{d-j}=0$ holds for $1\leq j \leq d-1$. Finally taking into
account that $(0,\ldots,0)$ belongs to
$\mathcal{L_I}\subset\Sigma_f\subset V_f$ we obtain $f_0=0$.
Therefore, we have the following result.
\begin{theorem}\label{teo:singular locus dim d-1 implies monomial}
  With notations as above, if the singular locus $\Sigma_f$ of $V_f$
  has dimension $d-1$, then $f_0=\cdots=f_{d-1}=0$ holds.
\end{theorem}

\subsection{The monomial case}
Fix a polynomial $f\in\fq[T]$ of degree $k+d<q-1$ with $k>d$ as in
(\ref{eq:polinomio f}) and consider the corresponding hypersurface
$V_f\subset\A^{k+1}$. Corollary \ref{coro:dimension singular locus
Vf} shows that the dimension of the singular locus of $V_f$ is at
most $d-1$. Furthermore, Theorem \ref{teo:singular locus dim d-1
implies monomial} asserts that, if the dimension of the singular
locus of $V_f$ is $d-1$, then the polynomial $f$ is necessarily the
monomial $f=T^{k+d}$. Our purpose in this section is to show that,
if the characteristic $p$ of $\fq$ satisfies the inequality $p>d+1$,
then this monomial does not generate a deep hole of the standard
Reed--Solomon code of dimension $k$ over $\fq$. This implies that,
for the sake of deciding the existence of deep holes, we may assume
without loss of generality that the singular locus of $V_f$ has
dimension at most $d-2$ when $p>d+1$ holds. As a first step in this
direction, we prove that, if the dimension of the singular locus of
$V_f$ is $d-1$, then $p$ divides $k+d$.
\begin{lemma}\label{lema:p divide a k+d}
Fix positive integers $k$ and $d$ with $k>d$. If the hypersurface
$V_d\subset\A^{k+1}$ associated to $T^{k+d}$ has a singular locus of
dimension $d-1$, then $p|(k+d)$.
\end{lemma}
\begin{proof}
We use the notations of the proof of Theorem \ref{teo:singular locus
dim d-1 implies monomial}. In the proof of Theorem \ref{teo:singular
locus dim d-1 implies monomial} we show that, if the singular locus
$\Sigma_d$ of $V_d$ has dimension $d-1$, then there exists a linear
variety
$$\mathcal{L_I}:=\mathrm{span}(\mathbf{v}^{(I_1)},\dots,
\mathbf{v}^{(I_{d-1})})$$
of dimension $d-1$ contained in $\Sigma_d$, where
$v_i^{I_j}\in\{0,1\}$ for $1\le i\le k+1$ and $1\le j\le d-1$, and
$\mathbf{v}^{(I_1)}+\dots+\mathbf{v}^{(I_{d-1})}=(1,\dots, 1)$. Let
$\lambda:=(\lambda_1,\dots,\lambda_{d-1})\in \A^{d-1}$ and let
$\mathbf{x}:=\sum_{j=1}^{d-1}\lambda_j \mathbf{v}^{(I_j)}$ be and
arbitrary point of $\mathcal{L_I}$. As in the proof of Theorem
\ref{teo:singular locus dim d-1 implies monomial}, we assume that
$x_i=\lambda_i$ $(1\leq i \leq d-1)$ and $\lambda_i\neq \lambda_j$
$(1\le i<j\le k+1)$ holds. Then the first $d-1$ equations of
(\ref{eq:singular locus de dimension d-1 primer sistema }) imply
that $\mathbf{0}\in\A^d$ is the unique solution of the following
nonsingular square system:
$$-B(\mathbf{x})\cdot\left(%
\begin{array}{c}
  f_1 \\
  \vdots \\
  f_{d-1} \\
\end{array}%
\right)=b(\mathbf{x})^{t},$$
where $B(\mathbf{x}):=\big((\partial H_j/\partial
X_i)(\mathbf{x}):{1\le i,j\le d-1}\big)$ and
$b(\mathbf{x}):=\big((\partial H_d/\partial X_j)(\mathbf{x}):1\le
j\le d-1\big)$.
In particular, the Cramer rule implies
\begin{equation}\label{eq:cramer rule proof monomial}
\det B^{(d-1)}(\mathbf{x})=0,
\end{equation}
where $B^{(d-1)}(\mathbf{x})\in \A^{(d-1)\times(d-1)}$ is the matrix
obtained by replacing the $(d-1)$th column of $B(\mathbf{x})$ by the
vector $b(\mathbf{x})$. We also recall that the matrix
$B^{(d-1)}(\mathbf{x})$ can be factored as in (\ref{eq:factorizacion
de Bj}), namely $B^{(d-1)}(\mathbf{x})=C(\mathbf{x})\cdot
H^{(d-1)}(\mathbf{x})$, where $C(\mathbf{x})$ is defined as in
(\ref{eq:matrix C}) and $H^{(d-1)}(\mathbf{x})\in\A^{d\times (d-1)}$
is the following matrix:
 $$
   H^{(d-1)}(\mathbf{x}):=\left(
    \begin{array}{ccccc}
     1 & H_1(\mathbf{x})\!\!\!\! & \cdots & \!\!\!\! H_{d-3}(\mathbf{x}) & H_{d-1}(\mathbf{x})
     \\[1ex]
      0 & 1& \cdots & \!\!\!\! H_{d-4}(\mathbf{x}) & H_{d-2}(\mathbf{x})
     \\
       & 0 &  \ddots & \vdots & \vdots
     \\
     \vdots  & \vdots & \ddots & 1 &H_2(\mathbf{x})
           \\[1ex]
      &  &  & 0 & H_1(\mathbf{x})
     \\[1ex]
     0 & 0 & \cdots & 0 & 1
    \end{array}
  \right).
  $$

We shall obtain an explicit expression of $\det
B^{(d-1)}(\mathbf{x})$ by applying the Cauchy--Binet formula to such
a factorization of $B^{(d-1)}(\mathbf{x})$. For this purpose, we
observe that $H^{(d-1)}(\mathbf{x})$ has  only two nonzero
$(d-1)\times (d-1)$ minors: the one corresponding to the submatrix
consisting of the first $d-1$ rows of $H^{(d-1)}(\mathbf{x})$, whose
value is equal to $H_1(\mathbf{x})$, and the one determined by the
rows $\{1, \ldots, d-2, d\}$ of $H^{(d-1)}(\mathbf{x})$, which is
equal to 1. Therefore, by the Cauchy-Binet formula we have
  $$
    \det B^{(d-1)}(\mathbf{x})= H_1 (\mathbf{x}) \cdot \det B(\mathbf{x})
   + \det \left(
    \begin{array}{ccccc}
     1 & x_1 & \cdots & x_{1}^{d-3} & x_{1}^{d-1}
     \\[1ex]
     1 & x_2 & \cdots & x_{2}^{d-3} & x_{2}^{d-1}
     \\ \vdots &\vdots&&\vdots&\vdots
     \\
     1 & x_{d-1} & \cdots & x_{d-1}^{d-3} & x_{d-1}^{d-1}
    \end{array}
  \right).
 $$
Combining (\ref{eq:cramer rule proof monomial}) with, e.g.,
\cite[Lemma 2.1]{Ernst00} or \cite[Exercise 280]{FaSo65}, we obtain
the following identity:
\begin{eqnarray*}
0&=&H_1(\mathbf{x})\cdot\det B(\mathbf{x})+(x_1+\cdots+ x_{d-1})\det
B(\mathbf{x})\\
&=&\det B(\mathbf{x})\cdot\big((\#I_1+1)\lambda_1+\cdots+
(\#I_{d-1}+1)\lambda_{d-1}\big).\end{eqnarray*}
Since $B(\mathbf{x})$ is a nonsingular matrix, we conclude that
\begin{equation}\label{eq:identity lambdas}
(\#I_1+1)\lambda_1+\cdots+ (\#I_{d-1}+1)\lambda_{d-1}=0
\end{equation}
holds for every $\lambda\in\A^{d-1}$ with $\lambda_i\not=\lambda_j$
for $i\not=j$, and hence for every $\lambda\in\A^{d-1}$. Substituting
$1$ for $\lambda_i$ in (\ref{eq:identity lambdas}), the statement of
the lemma follows.
\end{proof}
\begin{remark}
The conclusion in the statement of Lemma \ref{lema:p divide a k+d},
namely that $p|(k+d)$, is actually a rather weak consequence of
(\ref{eq:identity lambdas}). In addition to such a conclusion,
(\ref{eq:identity lambdas}) establishes strong restrictions on the
partitions $I$ of the linear varieties $\mathcal{L}_I$ contained in
the singular locus $\Sigma_f$ of a hypersurface $V_f$ with $\dim
\Sigma_f = d-1$. In particular, fix $i\in\{1,\dots,d-1\}$ and
substitute 1 for $\lambda_i$ and $0$ for any $\lambda_j$ with $j\neq
i$. Then (\ref{eq:identity lambdas}) implies $\# I_i\equiv -1\mod
p$.
\end{remark}

Now we are ready to prove the main result of this section, namely
that the assumption $p>d+1$ implies that any member $V_f$ of the
family of hypersurfaces which are relevant for the nonexistence of
deep holes has singular locus of dimension at most $d-2$.
\begin{proposition}\label{prop:monomial doesnt generate deep holes}
Let be given positive integers $k$ and $d$ with $k>d$, $p>d+1$, and
$q>k+d$. Assume further that $p|(k+d)$ holds. Let $\mathbf{w}_d\in
\fq^{q-1}$ be the word generated by the polynomial
$T^{k+d}\in\fq[T]$. Then $\mathbf{w}_d$ is not a deep hole of the
standard Reed--Solomon code $C$ of dimension $k$ over $\fq$.
\end{proposition}
\begin{proof} Write $q:=p^s$. The inequality $q>k+d\ge p$ implies $s>1$.
Consider the trace mapping $\mathrm{tr}_{\fq/\fp}:
\fq\rightarrow\fp$ defined by
$\mathrm{tr}_{\fq/\fp}(\alpha)=\sum_{i=0}^{s-1}\alpha^{p^{i}}$. It
is well--known that $\mathrm{tr}_{\fq/\fp}$ is a surjective
$\fp$--linear morphism. This in particular implies that there exist
$p^{s-1}$ elements in $\fq$ whose trace equals zero. Write $k+d=p\,
l$. Then the condition $q>k+d$ implies $p^{s-1}>l$, which in turn
shows that there exist $l$ pairwise--distinct elements
$b_1\ldots,b_l \in\fq^*$ with $\mathrm{tr}_{\fq/\fp}(b_i)=0.$

Since $\mathrm{tr}_{\fq/\fp}(b_i)=0$ holds for $1\leq i \leq l$, by
\cite[Theorem 3]{CoHe04} it follows that the Artin--Schreier
polynomial $g_{b_i}:=T^p-T-b_i\in\fq[T]$ has $p$ distinct roots in
$\fq^*$ for $1\leq i \leq l$. Furthermore, since $b_i\neq b_j$ holds
for $i\neq j$, we easily deduce that $g_{b_i}$ and $g_{b_j}$ have no
common roots. Therefore, the polynomial
\begin{equation}\label{eq:polynomial prod Artin-Schreier}
g:=\prod_{i=1}^lg_{b_i}=\prod_{i=1}^l(T^p-T-b_i)
\end{equation}
has $p\,l$ distinct roots in $\fq^*$. On the other hand,
$$g=T^{k+d}-lT^{p(l-1)+1}+\mathcal{O}(T^{p(l-1)})=T^{k+d}+h(T),$$
where $h:=lT^{p(l-1)+1}+\mathcal{O}(T^{p(l-1)})$ has degree at most
$p(l-1)+1$. Denote by $\mathbf{w}_h\in\fq^{q-1}$ the word generated
by the polynomial $h$. Since
$$p(l-1)+1=k+d-p+1\le k+d-(d+2)+1=k-1,$$
holds, we have that $\mathbf{w}_h$ is a codeword. The fact that the
polynomial $g$ of (\ref{eq:polynomial prod Artin-Schreier}) has
$p\,l>k$ distinct roots in $\fq^*$ implies ${\sf
d}(\mathbf{w}_d,\mathbf{w}_h)<q-1-k$ holds, where ${\sf d}$ denotes
the Hamming distance of $\fq^{q-1}$. We conclude that $\mathbf{w}_d$
is not a deep hole of the code $C$. This finishes the proof of the
proposition.
\end{proof}

\section{The main results}

We have shown that, if a given hypersurface $V_f$ has a
$q$--rational point with nonzero, pairwise--distinct coordinates,
then there are no deep holes of the standard Reed--Solomon code $C$
of dimension $k$ over $\fq$. Combining the results of Sections
\ref{section:geometry of Vf} and \ref{section:singular locus for
large characteristic}, we will obtain a lower bound for the number
of $q$--rational points of $V_f$ and an upper bound  for the number
of $q$--rational points of $V_f$ with a zero coordinate or at least
two equal coordinates. From these results we will establish a lower
bound for the number of $q$--rational points of $V_f$ as required.
This will allow us to obtain conditions on $q$, $d$ and $k$ which
imply the nonexistence of deep holes of the standard Reed--Solomon
code $C$.

As before, let be given positive integers $d$ and $k$ with $k>d$ and
$q-1>k+d$ and a polynomial $f:=T^{k+d}+f_{d-1}T^{k+d-1}+\cdots+
f_0T^k \in\fq[T]$. Consider the hypersurface $V_f\subset\A^{k+1}$
defined by the polynomial $H_f\in \fq[X_1,\ldots X_{k+1}]$
associated to $f$. According to Corollaries \ref{coro:dimension
singular locus Vf} and \ref{coro:Vf is abs irred}, the hypersurface
$V_f$ has a singular locus of dimension at most $d-1$ and is
absolutely irreducible.
%
%
\subsection{Estimates on the number of $q$--rational points of
hypersurfaces}
In what follows, we shall use an estimate on the number of
$q$--rational points of a projective $\fq$--hypersurface due to S.
Ghorpade and G. Lachaud (\cite{GhLa02a}; see also \cite{GhLa02}). In
\cite[Theorem 6.1]{GhLa02a} the authors prove that, for an
absolutely irreducible $\fq$--hypersurface $V\subset\Pp^{m+1}$ of
degree $d\ge 2$ and singular locus of dimension at most $s\ge 0$,
the number $\# V(\fq)$ of $q$--rational points of $V$ satisfies the
estimate
  \begin{equation}\label{eq:estimate rat points Ghorpade-Lachaud}
    |\# V(\fq)-p_{m}|\leq b_{m-s-1,d}\,q^{\frac{m+s+1}{2}}+C_{s,m}(V)q^{\frac{m+s}{2}},
  \end{equation}
where $p_{m}:=q^{m}+q^{m-1}+\cdots+ q+1$ is the cardinality of
$\Pp^{m}(\fq)$. Here  $b_{m-s-1,d}$ is the $(m-s-1)$th primitive
Betti number of any nonsingular hypersurface in $\Pp^{m-s}$ of
degree $d$, which is upper bounded by
  \begin{equation}\label{eq:upper bound betti number}
    b_{m-s-1,d} \le \frac{d-1}{d}\big((d-1)^{m-s}-(-1)^{m-s}\big)\le (d-1)^{m-s},
  \end{equation}
while  $C_{s,m}(V)$ is the sum
  $$
    C_{s,m}(V):=\sum_{i=m}^{m+s}b_{i,\ell}(V)+\varepsilon_i,
  $$
where $b_{i,\ell}(V)$ denotes the $i$th $\ell$--adic Betti number of
$V$ for a prime $\ell$ different from $p:=\mathrm{char}(\fq)$ and
$\varepsilon_i:=1$ for even $i$ and $\varepsilon_i:=0$ for odd $i$.
In \cite[Proposition 5.1]{GhLa02a} it is shown that
  \begin{equation}\label{eq:upper bound betti number Ghorpade Lachaud}
    C_{s,m}(V)\leq 18 (d + 3)^{m+2}.
  \end{equation}
This bound is a particular case of a bound for singular projective
complete intersections. Nevertheless, in our case it is possible to
slightly improve (\ref{eq:upper bound betti number Ghorpade
Lachaud}).
\begin{lemma}\label{lema:improved upper bound sums betti number}
If  $V\subset\Pp^{m+1}$  is  an absolutely irreducible hypersurface
of degree $d\ge 2$ and singular locus of dimension at most $s\ge 0$,
then we have the following bound:
  \begin{equation}\label{eq:upper bound sums betti numbers}
     C_{s,m}(V) \le  6(d+2)^{m+2}
  \end{equation}
\end{lemma}
  \begin{proof}
Let $E(n,d)$ be an upper bound for the Euler characteristic of the
affine part and the part at infinity of $V$.  Let $A(n,d)$ be the
number
       $$
         A(n, d):= E(n,d) + 2 + 2 \sum_{j=1}^{n-1}E(j,d).
       $$
Then the Katz inequality \cite[Theorem 3]{Katz01} asserts that
       $$
         C_{m,s}(V) \leq 1 + \sum_{n=1}^{m+1}  \big( 1 +  A(n+1,d+1) \big).
       $$

As a consequence of \cite[Theorem 5.27]{AdSp88}, we have
       $$
         E(n,d):=2 (d+1)^n.
       $$
These results and elementary calculations easily imply
(\ref{eq:upper bound sums betti numbers}).
  \end{proof}

Combining (\ref{eq:estimate rat points Ghorpade-Lachaud}) with
(\ref{eq:upper bound betti number}) and Lemma \ref{lema:improved
upper bound sums betti number} we obtain an explicit upper bound for
the number of $q$--rational points of singular projective
$\fq$--hypersurfaces. More precisely, if $V\subset\Pp^{m+1}$ is an
absolutely irreducible $\fq$--hypersurface of degree $d\ge 2$ and
singular locus of dimension at most $s\ge 0$, then the number of
$q$--rational points of $V$ satisfies the estimate
  \begin{equation}\label{eq:estimate rat points Ghorpade-Lachaud explicita}
    |\# V(\fq)-p_{m}|\leq (d-1)^{m-s}\,q^{\frac{m+s+1}{2}}+ 6 (d+2)^{m+2} q^{\frac{m+s}{2}}.
  \end{equation}

The first step towards our main result is to obtain a lower bound on
the number of $q$--rational points of the hypersurface $V_f$. For
this purpose, combining Corollary \ref{coro:Vf is abs irred} and
\cite[Chapter I, Proposition 5.17]{Kunz85} we conclude that the
projective closure $\mathrm{pcl}(V_f)\subset\Pp^{k+1} $ of $V_f$ is
an absolutely irreducible hypersurface which is defined over $\fq$.
Furthermore, from Corollary \ref{coro:dimension singular locus Vf}
and Proposition \ref{prop:dimension singularidades en el infinito}
we deduce that the singular locus of $\mathrm{pcl}(V_f)$ has
dimension at most $d-1$. Therefore from (\ref{eq:estimate rat points
Ghorpade-Lachaud explicita}) we deduce the following estimate:
  \begin{equation}\label{eq:estimate rat points proj clos}
    |\#\mathrm{pcl}(V_f)(\fq)-p_k|\leq (d-1)^{k-d+1}\,q^{\frac{k+d}{2}}+6(d+2)^{k+2}q^{\frac{k+d-1}{2}}.
  \end{equation}

Our next result provides a lower bound on the number of
$q$--rational zeros of the affine hypersurface $V_f$.
\begin{proposition}\label{prop:lower bound rat points Vf}
  Let be given positive integers $d$ and $k$ with $k>d\ge 2$ and $q-1>k+d$. Then
  the number of $q$--rational points of the hypersurface $V_f$
  satisfies the following inequality:
   $$
     \# V_f(\fq) \ge  q^k-2(d-1)^{k-d+1}q^{\frac{k+d}{2}} - 7(d+2)^{k+2}q^{\frac{k+d-1}{2}}.
   $$
\end{proposition}
 \begin{proof}
   Since we are interested in the $q$--rational points of $V_f$, we
   discard the points of $\mathrm{pcl}(V_f)(\fq)$ lying in the
   hyperplane at infinity $\{X_0=0\}$. Since $\mathrm{pcl}(V_f)$ is the
   locus of zeros of the polynomial $H_f^h=H_d+f_{d-1}H_{d-1}X_0+\cdots+
   f_0X_0^d\in\fq[X_0,\ldots, X_{k+1}]$, we conclude
   $$
    \mathrm{pcl}(V_f)(\fq)\cap\{X_0=0\}=\{\mathbf{x}\in \Pp^k(\fq):H_d(\mathbf{x})=0\}.
   $$

According to Proposition \ref{prop:dimension singularidades en el
infinito}, the projective $\fq$--hypersurface
$V_f^\infty\subset\Pp^k$ defined by $H_d$ has a singular locus of
dimension at most $d-2$. Applying (\ref{eq:estimate rat points
Ghorpade-Lachaud explicita}) we obtain
  \begin{equation}\label{eq:estimate rat points at infinity}
     |\#V_f^\infty(\fq)-p_{k-1}|\leq (d-1)^{k-d+1}\, q^{\frac{k+d-2}{2}}+6(d+2)^{k+1}\,q^{\frac{k+d-3}{2}}.
  \end{equation}
Combining (\ref{eq:estimate rat points proj clos}) and (\ref{eq:estimate rat points at infinity}) we have:
  \begin{eqnarray*}
      \# V_f(\fq)-q^k\!&\!\!\!\!\! = &\!\!\!\!
      \big(\#\mathrm{pcl}(V_f)(\fq)-p_k\big)-\big(\#V_f^\infty(\fq)-p_{k-1}\big)
      \\[1ex]
        &\!\!\!\!\! \ge  &\!\!\!\! -(d-1)^{k-d+1} q^{\frac{k+d}{2}} (1\!+\!q^{-1})-
        6 (d+2)^{k+2} q^{\frac{k+d-1}{2}}\big(1\!+\!(q(d\!+\!2))^{-1}\big).
  \end{eqnarray*}
From this lower bound the inequality of the statement easily
follows.
\end{proof}

Next we obtain an upper bound on the number of $q$--rational points
of the hypersurface $V_f$ which are not useful in connection with
the existence of deep holes, namely those with a zero coordinate or
at least two equal coordinates. We begin with the case of the points
with a zero coordinate.
\begin{proposition} \label{prop:zeros with a zero coordinate}
With hypotheses as in Proposition \ref{prop:lower bound rat points
Vf}, the number $N_1$ of $q$--rational points of $V_f$ with a zero
coordinate satisfies the following inequality:
  $$
    N_1 \le (k+1)\left(q^{k-1}+2(d-1)^{k-d}q^{\frac{k+d-1}{2}}
           +7(d+2)^{k+1}q^{\frac{k+d-2}{2}}\right).
  $$
\end{proposition}

  \begin{proof}
        Let $\mathbf{x}:=(x_1,\ldots, x_{k+1})$ be a point of $V_f$ with a zero coordinate. Without loss
        of generality we may assume $x_{k+1}=0$. Hence,
$\mathbf{x}$ is a $q$--rational point of the intersection
$W_{k+1}:=V_f\cap \{X_{k+1}=0\}$. Observe that $W_{k+1}$
        is the $\fq$--hypersurface of the linear space $\{X_{k+1}=0\}$ defined by the polynomial
        $G_f(\Pi_{1}^{k},\ldots,\Pi_{d}^{k})$, where $\Pi_{1}^{k},\ldots,\Pi_{d}^{k}$ are
        the first $d$ elementary symmetric polynomials of the ring
        $\fq[X_1,\ldots, X_{k}]$. Then Theorem
        \ref{teo:caracterizacion de los puntos singulares} shows that $W_{k+1}$
        has a singular locus of dimension at most $d-1$. Furthermore,
        Proposition \ref{prop:dimension singularidades en el infinito}
        implies that the singular locus of $W_{k+1}$ at infinity has also
        dimension at most $d-2$. As a consequence, arguing as in the proof
        of Proposition \ref{prop:lower bound rat points Vf} we
        obtain
        \begin{eqnarray*}
          \#W_{k+1}(\fq)-q^{k-1}\!\!\!\!&=\!\!\!\!& \big(\#\mathrm{pcl}(W_{k+1})(\fq)-p_{k-1}\big)
            -\big(\#W_{k+1}^\infty(\fq)-p_{k-2}\big)
          \\[1ex]
             &\le\!\!\!\!& (d-1)^{k-d}q^{\frac{k+d-1}{2}}\!  + 6 (d+2)^{k+1}q^{\frac{k+d-2}{2}}
          \\[1ex]
             &  &+(d-1)^{k-d}q^{\frac{k+d-3}{2}}\!  + 6 (d+2)^{k}q^{\frac{k+d-4}{2}}.
        \end{eqnarray*}
        %
        %
        Therefore, we have the upper bound
        \begin{equation}\label{eq:upper bound zeros with a zero coordinate}
           \#W_{k+1}(\fq) \le  q^{k-1}+2(d-1)^{k-d}q^{\frac{k+d-1}{2}}
                               + 7(d+2)^{k+1}q^{\frac{k+d-2}{2}}.
        \end{equation}
Adding the upper bounds of the $q$--rational points of the varieties
$W_i:=V_f\cap\{X_i=0\}$ for $1\le i\le k+1$, the proposition
follows.
  \end{proof}

Next we consider the number of $q$--rational points of $V_f$ with
two equal coordinates.
\begin{proposition} \label{prop:zeros with two equal coordinates}
With hypotheses as in Proposition \ref{prop:lower bound rat points
Vf}, the number $N_2$ of $q$--rational points of $V_f$ with at least two equal
coordinates satisfies the following inequality:
  $$
    N_2 \leq  \frac{(k+1)k}{2}\Big(q^{k-1}+2(d-1)^{k-d}q^{\frac{k+d-1}{2}}
    +7(d+2)^{k+1}q^{\frac{k+d-2}{2}}\Big).
  $$
\end{proposition}
  \begin{proof}
Let $\mathbf{x}:=(x_1,\ldots, x_{k+1})\in V_f(\fq)$ be a point
having two distinct coordinates with the same value.
    Without loss of generality
    we may assume that $x_k = x_{k+1}$ holds. 
    Then $\mathbf{x}$ is a $q$--rational point of the hypersurface
    $W_{k,k+1}\subset \{X_k=X_{k+1}\}$ defined by the polynomial
    $G_f(\Pi_{1}^*,\ldots,\Pi_{d}^*)\in \fq[X_1,\ldots,X_k]$, where
    $\Pi_i^*:=\Pi_{i}(X_1,\ldots, X_k,X_k)$ is the polynomial of $\fq[X_1,\ldots, X_k]$
    obtained by substituting $X_k$ for $X_{k+1}$ in the $i$th elementary
    symmetric polynomial of $\fq[X_1,\ldots, X_{k+1}]$. Observe that
\begin{equation}\label{eq:pi_i en x_k igual x_(k+1)}
      \Pi_{i}^* = \Pi_i^{k-1} + 2X_{k}\cdot \Pi_{i-1}^{k-1} + X_k^2 \cdot \Pi_{i-2}^{k-1}
\end{equation}
    where $\Pi_j^l$ denotes the $j$th elementary symmetric polynomial of
    $\fq[X_1,\ldots, X_l]$ for $1\le j\le d$ and $1\le l\le k+1$.

    We claim that the singular locus of
    $\mathrm{pcl}(W_{k,k+1})$ and the singular locus of $W_{k,k+1}$ at
    infinity have dimension at most $d-1$ and $d-2$, respectively.
    First, if the characteristic  $p$ of $\fq$ is greater than $2$,
    then using (\ref{eq:pi_i en x_k igual x_(k+1)}) it can be proved
    that all the maximal minors of the Jacobian matrix $(\partial\Pi^*_i/\partial
    X_j)_{1\le i\le d,1\le j\le k}$ are equal to the corresponding
    minors of the Jacobian matrix $(\partial\Pi^k_i/\partial
    X_j)_{1\le i\le d,1\le j\le k}$. Then the proofs of Theorem
    \ref{teo:caracterizacion de los puntos singulares} and Proposition
    \ref{prop:dimension singularidades en el infinito} go through with
    minor corrections and show our claim.

    For $p=2$, from (\ref{eq:pi_i en x_k igual x_(k+1)}) we see that
    the first partial derivative of $\Pi_j^*$ with respect to $X_k$ is
    equal to zero. Thus, each nonzero maximal minor of the
    Jacobian matrix $(\partial{\Pi_{i}^*}/\partial{X_j})_{1\le i\le d,1\le j\le k}$
    is a Vandermonde determinant depending on $d$ of the indeterminates
    $X_1, \ldots, X_{k-1}$. In particular, the vanishing of all
    these minors does not impose any condition on the variable
    $X_k$.

    Let $\Sigma_{k,k+1}$ denote the singular locus of $W_{k,k+1}$.
    Arguing as in the proof of Theorem \ref{teo:caracterizacion de los puntos singulares},
    we have the following inclusion (see Remark \ref{obs: caracterizacion de las variedades lineales}):
\begin{equation}\label{eq: singular locus sigma k k+1}
      \Sigma_{k,k+1}\subset \bigcup_{\mathcal{I}}\mathcal{L_{I}},
\end{equation}
where $\mathcal{I}:=\{I_1,\ldots,I_{d-1}, I_d\}$ runs over all the
partitions of $\{1,\ldots,k+1\}$ into $d$ nonempty subsets
$I_j\subset \{1,\ldots,k+1\}$ such that $I_j\subset\{1,\dots,k-1\}$
for $1\le j\le d-1$ and $I_d:=\{k,k+1\}$, $\mathcal{L_{I}}$ is the
linear variety
    $$
      \mathcal{L_{I}}:= \mathrm{span}(\mathbf{v}^{(I_1)},\ldots,\mathbf{v}^{(I_{d})})
    $$
spanned by the vectors $\mathbf{v}^{(I_j)}:=
(v_1^{(I_j)},\ldots,v_{k+1}^{(I_j)})$ defined by $v_m^{(I_j)}:=1$
for $m\in I_j$ and $v_m^{(I_j)}:=0$ for $m\notin I_j$. It follows
that, if $\Sigma_{k,k+1}$ has dimension $d$, then it contains a
linear variety $\mathcal{L_{I}}$ as above.

Suppose that $\Sigma_{k,k+1}$ has dimension $d$. Following the proof
of Theorem \ref{teo:singular locus dim d-1 implies monomial} we
conclude that $f$ is the monomial $T^{k+d}$, and thus $H_f=H_d$
holds. Fix $\mathcal{I}:=\{I_1,\ldots,I_{d-1}, I_d\}$ as above and
consider the corresponding $d$--dimensional linear variety
$\mathcal{L_{I}}$. We claim that $\mathcal{L_{I}}$ intersects
$\Sigma_{k,k+1}$ properly. Observe that, combining this claim with
(\ref{eq: singular locus sigma k k+1}), we easily deduce that
$\dim\Sigma_{k,k+1}\le d-1$, since each variety $\mathcal{L_{I}}$ is
absolutely irreducible and each irreducible component of
$\Sigma_{k,k+1}$ is a proper subvariety of a suitable
$\mathcal{L_{I}}$.

In order to prove our claim, consider the line
$\ell_\lambda:=\{\mathbf{v}_\lambda:=
(0,\dots,0,\lambda,\lambda)\in\A^{k+1}:\lambda\in\A^1\}
\subset\mathcal{L_I}$. Observe that
$\ell_\lambda\cap\Sigma_{k,k+1}=\{\mathbf{v}_\lambda\in\A^{k+1}:H_d(
\mathbf{v}_\lambda)=0,\nabla H_d(\mathbf{v}_\lambda)=\mathbf{0}\}$.
From the identities $\Pi_j(\mathbf{v}_\lambda)=0$ $(j\notin\{0,2\})$
and $\Pi_2(\mathbf{v}_\lambda)=\lambda^2$ and Proposition
\ref{Proposicion: formula explicita para la hipersuperficie H_0d},
we conclude that $H_d(\mathbf{v}_\lambda)=\pm\lambda^d$ for even $d$
and $H_{d-1}(\mathbf{v}_\lambda)=\pm\lambda^{d-1}$ for odd $d$.
Furthermore, from Lemma \ref{Lema: derivadas de los H0d} we obtain
$(\partial H_d/\partial
X_1)(\mathbf{v}_\lambda)=H_{d-1}(\mathbf{v}_\lambda)=\pm\lambda^{d-1}$
for odd $d$. In both cases, the identities
$H_d(\mathbf{v}_\lambda)=(\partial H_d/\partial
X_1)(\mathbf{v}_\lambda)=0$ imply $\lambda=0$. This shows that
$\ell_\lambda\subset\mathcal{L_I}$ intersects properly
$\Sigma_{k,k+1}$ and shows our claim.

Finally, arguing as in Proposition \ref{prop:dimension
singularidades en el infinito} we conclude that the singular locus
of $W_{k,k+1}$ at infinity has dimension at most $d-2$.

Summarizing, we have that, independently of the characteristic $p$
of $\fq$, the singular locus of $\mathrm{pcl}(W_{k,k+1})$ and the
singular locus of $W_{k,k+1}$ at infinity have dimension at most
$d-1$ and $d-2$. Then, following the proof of Proposition
    \ref{prop:zeros with a zero coordinate} we obtain:
    \begin{equation}\label{eq:upper bound zeros with two equal coordinates}
    \#W_{k,k+1}(\fq)\le  q^{k-1}+2(d-1)^{k-d}q^{\frac{k+d-1}{2}}
       +7(d+2)^{k+1}q^{\frac{k+d-2}{2}}.
    \end{equation}
    From  (\ref{eq:upper bound zeros with two equal coordinates})
    we deduce the statement of the proposition.
 \end{proof}
%
%
\subsection{Results of nonexistence of deep holes}
Now we are ready to prove the main results of this paper. Fix $q$,
$k$ and $d\ge 3$ with $q>k+d$ and consider the standard
Reed--Solomon code $C$ of dimension $k$ over $\fq$. From Section
\ref{section:introduction} we have that a polynomial
$f:=T^{k+d}+f_{d-1}T^{k+d-1}+\cdots+f_{0}T^k$ does not generate a
deep hole of the code $C$ if and only if the corresponding
hypersurface $V_f\subset\A^{k+1}$ has a $q$--rational point with
nonzero, pairwise--distinct coordinates. Combining Propositions
\ref{prop:lower bound rat points Vf}, \ref{prop:zeros with a zero
coordinate} and \ref{prop:zeros with two equal coordinates} we
conclude that the number $N$ of such points satisfies the following
inequality:
\begin{equation}\label{eq:lower bound useful points}
  \begin{split}
  N\ge  q^k-\dfrac{(k+1)(k+2)}{2}q^{k-1}-2(d-1)^{k-d}q^{\frac{k+d}{2}}\bigg(d-1+
    \dfrac{(k+1)(k+2)}{2q^{\frac{1}{2}}}\bigg) \\
    -7(d+2)^{k+1}q^{\frac{k+d-1}{2}}\bigg(d+2+
    \dfrac{(k+1)(k+2)}{2q^{\frac{1}{2}}}\bigg).
  \end{split}
\end{equation}
Therefore, the polynomial $f$ does not generate a deep hole of the
code $C$ if the right--hand side of (\ref{eq:lower bound useful
points}) is a positive number.

Suppose that $q$, $k$ and $d\ge 3$ satisfy the following conditions:
 \begin{equation}\label{eq:first condition on q}
   q>(k+1)^2,\quad k>3d.
 \end{equation}
Since $k\geq 10$,  it follows that
  $
    \frac{3}{4}(k+1)(k+2)\leq (k+1)^{2}<q
  $
holds. Therefore, we have
  $
    q-\frac{1}{2}(k+1)(k+2)>{q}/{3},
  $
which implies
  $$
    q^k-\dfrac{(k+1)(k+2)}{2}q^{k-1}=q^{k-1}\bigg(q-\dfrac{(k+1)(k+2)}{2}\bigg)> \frac{q^k}{3}.
  $$
Hence, the right--hand side of (\ref{eq:lower bound useful points})
is positive if the following condition holds:
 \begin{equation}
 \label{eq: first inequality of the proof of lower bound useful points}
   \begin{array}{ccl}
    \dfrac{q^k}{3}  & \geq  &   2(d-1)^{k-d}q^{\frac{k+d}{2}}\bigg(d-1+
    \dfrac{(k+1)(k+2)}{2q^{\frac{1}{2}}}\bigg)\qquad \\ &  &\qquad+\, 7(d+2)^{k+1}q^{\frac{k+d-1}{2}}\bigg(d+2+
    \dfrac{(k+1)(k+2)}{2q^{\frac{1}{2}}}\bigg).
   \end{array}
 \end{equation}
Taking into account that  $k+1<q^{\frac{1}{2}}$, we conclude that
(\ref{eq: first inequality of the proof of lower bound useful
points}) can be replaced by the following condition:
  $$
    \begin{array}{c}
      \dfrac{q^k}{3}\ge 2(d-1)^{k-d}q^{\frac{k+d}{2}}\left(d-1+\frac{k+2}{2}\right)
      +  7(d+2)^{k+1}q^{\frac{k+d-1}{2}}\left(d+2+\frac{k+2}{2}\right).
    \end{array}
 $$
From $d\leq\frac{k-1}{3}$ we obtain
  $
    d+2+\frac{k+2}{2}\leq k+1,
  $
and therefore we conclude that the right--hand side of
(\ref{eq:lower bound useful points}) is positive if
  $$
   \dfrac{q^k}{3}\ge 2(d-1)^{k-d}(k-2)q^{\frac{k+d}{2}}  +7(d+2)^{k+1}(k+1)q^{\frac{k+d-1}{2}},
  $$
or equivalently if
  \begin{equation}
     q^k\ge 6(d-1)^{k-d}(k-2)q^{\frac{k+d}{2}}
     +21(d+2)^{k+1}(k+1)q^{\frac{k+d-1}{2}},
  \end{equation}
holds.
Furthermore, this condition is in turn implied by the following
conditions:
$$    \dfrac{q^k}{8}\ge 6(d-1)^{k-d}(k-2)q^{\frac{k+d}{2}},\quad
\dfrac{7 q^k}{8}\ge
    21(d+2)^{k+1}(k+1)q^{\frac{k+d-1}{2}},$$
which can be rewritten as
  \begin{equation}\label{eq:first condition on q bis}
     q^k\ge 48(d-1)^{k-d}(k-2)q^{\frac{k+d}{2}},\quad q^k\ge 24(d+2)^{k+1}(k+1)q^{\frac{k+d-1}{2}}.
  \end{equation}
The first inequality is equivalent to the following inequality:
  $$
    q\ge (48(k-2))^{\frac{2}{k-d}}(d-1)^2.
  $$
From (\ref{eq:first condition on q}) one easily concludes that
$3(k-d)\ge 2k+1$ holds. Since the function $k\mapsto
\big(48(k-2)\big){}^{6/(2k+1)}$ is decreasing, taking into account
that $k\ge 10$ holds we deduce that a sufficient condition for the
fulfillment of the inequality above is
\begin{equation}\label{eq:first part first condition}
     q>6d^2.
\end{equation}
Next we consider the second inequality of (\ref{eq:first condition
on q bis}). First, we observe that this inequality can be expressed
as follows:
\begin{equation}\label{eq:second part first condition inexplicit}
    q>(24(k+1))^{\frac{2}{k-d+1}}\Big(\frac{d+2}{d}\Big)^{2+\frac{2d}{k-d+1}}d^{2+\frac{2d}{k-d+1}}
\end{equation}
From (\ref{eq:first condition on q}) we deduce $3(k-d+1)\ge 2k+4$.
Taking into account that the function $k\mapsto
\big(24(k+1)\big){}^{3/(k+2)}$ is decreasing, in particular for $k \geq 12$ (and thus for $d\ge 4$),
we see that (\ref{eq:second part first condition
inexplicit}) is satisfied if the following condition holds:
  \begin{equation}\label{eq:second part first condition}
    q>14\,d^{2+2d/(k-d)}.
  \end{equation}

Combining (\ref{eq:first condition on q}), (\ref{eq:first part first
condition}) and (\ref{eq:second part first condition}) we conclude
that (\ref{eq:first condition on q}) and (\ref{eq:second part first
condition}) yield a sufficient condition for the nonexistence of
deep holes. Finally, starting from (\ref{eq:lower bound useful
points}) one easily sees that (\ref{eq:first condition on q}) and
(\ref{eq:second part first condition}) yield a sufficient condition
for the nonexistence of deep holes for $d=3$.

\begin{theorem}\label{teo:principal - version precisa}
Let  $k$ and $d$ be integers with $k>d\ge 3$ and $q-1>k+d$, and let
${C}$ be the standard Reed--Solomon code of dimension $k$ over
$\fq$. Let be given a real number $\epsilon$ with $0<\epsilon<1$ and
let $\mathbf{w}$ be a word generated by a polynomial $f\in\fq[T]$ of
degree $k+d$. If the conditions
  $$
    q>\max\{(k+1)^2,14\,d^{2+\epsilon}\},\quad k\ge d\Big(\frac{\,2}{\epsilon}+1\Big)
  $$
hold, then $\mathbf{w}$ is not a deep hole of ${C}$.
\end{theorem}

We remark that in \cite{LiWa08b} it is shown that, for $d=1$, $k>2$
and $q>k+3$, polynomials of degree $k+1$ do not generate deep holes
of the standard Reed--Solomon code $C$. On the other hand, a similar
result as in Theorem \ref{teo:principal - version precisa} can be
obtained for $d=2$ with our approach, namely that for a suitable
constant $M_1>14$, if the conditions
$q>\max\{(k+1)^2,M_1\,2^{2+\epsilon}\}$ and $k\ge 2(2/{\epsilon}+1)$
hold, then no polynomial of degree $k+2$ generates a deep hole of
$C$. Finally, since the dimension of the singular locus of
$\mathrm{pcl}(V_f)$ is at most $d-1 \leq k-2$, by the Serre
criterion of normality (see, e.g., \cite[Theorem 18.15]{Eisenbud95})
we have that $\mathrm{pcl}(V_f)$ is a normal variety. Hence,
combining Theorem \ref{teo:caracterizacion de los puntos singulares}
and Proposition \ref{prop:dimension singularidades en el infinito}
with the estimates on the number of $q$--rational points of a normal
complete intersection defined over $\fq$ of \cite{CaMa07}, it can be
proved that for $k>d$, $q-1>k+d$ and $q>\max\{(k+1)^2,M_2\,d^4\}$
for a suitable constant $M_2$, no polynomial of degree $k+d$
generates a deep hole of $C$.
%
%
\subsection{Nonexistence of deep holes for $\mathrm{char}(\fq)>d+1$}
Finally, we briefly indicate what we obtain under the assumption
that the characteristic $p$ of $\fq$ satisfies the inequality
$p>d+1$. Fix $q$, $k$ and $d\ge 3$ with $q-1>k+d$, $k>d$ and $p>d+1$
and consider the standard Reed--Solomon code $C$ of dimension $k$
over $\fq$.

Fix $f:=T^{k+d}+f_{d-1}T^{k+d-1}+\cdots+f_0T^k\in\fq[T]$. First
suppose that the singular locus of the hypersurface $V_f$ associated
to $f$ has dimension $d-1$. By Theorem \ref{teo:singular locus dim
d-1 implies monomial} we have that $f$ is the monomial $T^{k+d}$.
Furthermore, from Lemma \ref{lema:p divide a k+d} it follows that
$p|(k+d)$. Then Proposition \ref{prop:monomial doesnt generate deep
holes} shows that the monomial $T^{k+d}$ does not generate a deep
hole of $C$. Therefore, we may assume without loss of generality
that $V_f$ has a singular locus of dimension at most $d-2$. As a
consequence, arguing as in the proofs of Propositions
\ref{prop:lower bound rat points Vf}, \ref{prop:zeros with a zero
coordinate} and \ref{prop:zeros with two equal coordinates} we
obtain the following bounds:
\begin{eqnarray*}
  \# V_f(\fq) & \ge & q^k-2(d-1)^{k-d+2}q^{\frac{k+d-1}{2}}-7(d+2)^{k+2}q^{\frac{k+d-2}{2}},
  \\[1ex]
  N_1&\le & \frac{(k+1)(k+2)}{2}\Big(q^{k-1}+2(d-1)^{k-d+1}q^{\frac{k+d-2}{2}}
   +7(d+2)^{k+1}q^{\frac{k+d-3}{2}}\Big),
\end{eqnarray*}
where $N_1$ denotes the number of $q$--rational points of $V_f$
having a zero coordinate or at least two equal coordinates. Hence we
have that the number $N$ of $q$--rational points of $V_f$ with
nonzero, pairwise--distinct coordinates satisfies the following
inequality:
 \begin{equation}\label{eq:lower bound useful points char greater than d}
  \begin{split}
 N\ge
q^k-\dfrac{(k\!+\!1)(k\!+\!2)}{2}q^{k-1}-2(d-1)^{k-d+1}q^{\frac{k+d-1}{2}}\bigg(d-1+
    \dfrac{(k\!+\!1)(k+2)}{2q^{\frac{1}{2}}}\bigg) \\
    -7(d+2)^{k+1}q^{\frac{k+d-2}{2}}\bigg(d+2+
    \dfrac{(k\!+\!1)(k+2)}{2q^{\frac{1}{2}}}\bigg).
  \end{split}
 \end{equation}
Suppose that $q$, $k$ and $d\ge 4$ satisfy the following conditions:
 \begin{equation}\label{eq:first condition on q char greater than d}
   q>(k+1)^2,\quad k>3(d-1).
 \end{equation}
Then the right--hand side of (\ref{eq:lower bound useful points char
greater than d}) is positive if
  \begin{equation}\label{eq:first condition on q bis char greater than d}
   q^k\ge\max\big\{48(d-1)^{k-d+1}(k+1)q^{\frac{k+d-1}{2}},\,24(d+2)^{k+1}(k+4)q^{\frac{k+d-2}{2}}\big\}.
  \end{equation}
With similar arguments as in the proof of Theorem \ref{teo:principal
- version precisa} we conclude that (\ref{eq:first condition on q
bis char greater than d}) is satisfied if the following condition
holds:
\begin{equation}\label{eq:second part first condition char greater than d}
    q>20\,d^{2+2d/(k-d)}.
\end{equation}
On the other hand, starting from (\ref{eq:lower bound useful points
char greater than d}) one easily sees that (\ref{eq:second part
first condition char greater than d}) yields a sufficient condition
for the nonexistence of deep holes for $d=3$. Summarizing, we have
the following result.
\begin{theorem}\label{teo:principal char greater than d}
Let  $k$ and $d$ be integers with $k>d\ge 3$ and $q-1>k+d$, and let
${C}$ be the standard Reed--Solomon code of dimension $k$ over
$\fq$. Let be given a real number $\epsilon$ with $0<\epsilon<1$ and
let $\mathbf{w}$ be a word generated by a polynomial $f\in\fq[T]$ of
degree $k+d$. If $\mathrm{char}(\fq)>d+1$ and the conditions
$$q>\max\{(k+1)^2,20\,d^{2+\epsilon}\},\quad k\ge (d-1)\Big(\frac{\,2}{\epsilon}+1\Big)$$
hold, then $\mathbf{w}$ is not a deep hole of ${C}$.
\end{theorem}

\end{document}